\documentclass[a4paper,11pt]{article}

\usepackage[english]{babel}
\usepackage[latin1]{inputenc}

\usepackage{lmodern}
\usepackage{textcomp}
\usepackage{enumerate}
\usepackage{graphicx,epsfig,psfrag}
\usepackage{mathrsfs}
\usepackage{amsmath,amsfonts,amssymb}
\usepackage{color}
\usepackage{xcolor}
\usepackage{fancybox}
\usepackage{pgf,tikz}
\usetikzlibrary{arrows} 
\usepackage{animate}
\definecolor{rltred}{rgb}{0.75,0,0} 
	\definecolor{rltgreen}{rgb}{0,0.5,0}
	\definecolor{oneblue}{rgb}{0,0,0.75}
	\definecolor{marron}{rgb}{0.64,0.16,0.16}
	\definecolor{forestgreen}{rgb}{0.13,0.54,0.13}
	\definecolor{purple}{rgb}{0.62,0.12,0.94}
	\definecolor{dockerblue}{rgb}{0.11,0.56,0.98}
	\definecolor{freeblue}{rgb}{0.25,0.41,0.88}
	\definecolor{myblue}{rgb}{0,0.2,0.4}
	\definecolor{ccqqtt}{rgb}{0.8,0,0.2}

\newtheorem{lem}{Lemma}[section]

\newtheorem{prop}{Proposition}[section]
\newtheorem{defin}{Definition}[section]
\newtheorem{theo}{Theorem}[section]

\newenvironment{dem}{\vskip 2mm\noindent \textbf{Proof:}}{\hfill $\square$ \vskip 2mm \noindent}

\title{Absence of percolation for Poisson outdegree-one graphs}

\newcommand{\1}{\ensuremath{\mbox{\rm 1\kern-0.23em I}}}
\newcommand{\C}{\mathscr{C}}
\newcommand{\Rd}{\mathbb{R}^{d}}
\newcommand{\Rde}{\mathbb{R}^{2}}
\newcommand{\Z}{\mathbb{Z}}

\renewcommand{\H}{\text{Hex}}

\begin{document}
\maketitle
\begin{center}
David Coupier, David Dereudre and Simon Le Stum\\
Laboratoire paul Painlev\'e, University Lille 1
\end{center}

\begin{center}
\textbf{Abstract}
\end{center}
A Poisson outdegree-one graph is an oriented graph based on a Poisson point process such that each vertex has only one outgoing edge. The paper focuses on the absence of percolation for such graphs. Our main result is based on two assumptions. The Shield assumption ensures that the graph is locally determined with possible random horizons. The Loop assumption ensures that any forward branch of the graph merges on a loop provided that the Poisson point process is augmented with a finite collection of well-chosen points. Several models satisfy these general assumptions and inherit in consequence the absence of percolation. In particular, we solve in Theorem \ref{mainTH} a conjecture by Daley et al. on the absence of percolation for the line-segment model (Conjecture 7.1 of \cite{daley2014two}, discussed in \cite{hirsch2016absence} as well). In this planar model, a segment is growing from any point of the Poisson process and stops its growth whenever it hits another segment. The random directions are picked independently and uniformly on the unit sphere. Another model of geometric navigation is presented and also fulfills the Shield and Loop assumptions.

\section{Introduction}
\label{SectionIntro}

Consider the classical nearest neighbour graph based on a planar homogeneous Poisson point process in which each point is simply connected to its nearest neighbour. The absence of percolation for this model is due to the fact that almost surely a homogeneous Poisson point process contains no descending chain. By a descending chain, we mean an infinite sequence $x_{1},x_{2},...$ of points of the process for which $|x_{i-1}-x_{i}|\ge |x_{i}-x_{i+1}| $ for all $i\ge 2$.  Daley and Last have shown in \cite{daley2005descending} that the absence of percolation for the lilypond model can also be obtained as a consequence of the descending chain argument. In this model a ball is growing with unit rate from any Poisson point and stops its growth when it hits another ball. Note that the finite cluster property for this model has first been proved in \cite{haggstrom1996nearest}.

\begin{figure}[!h]
\begin{center}
\includegraphics[width=6.5cm]{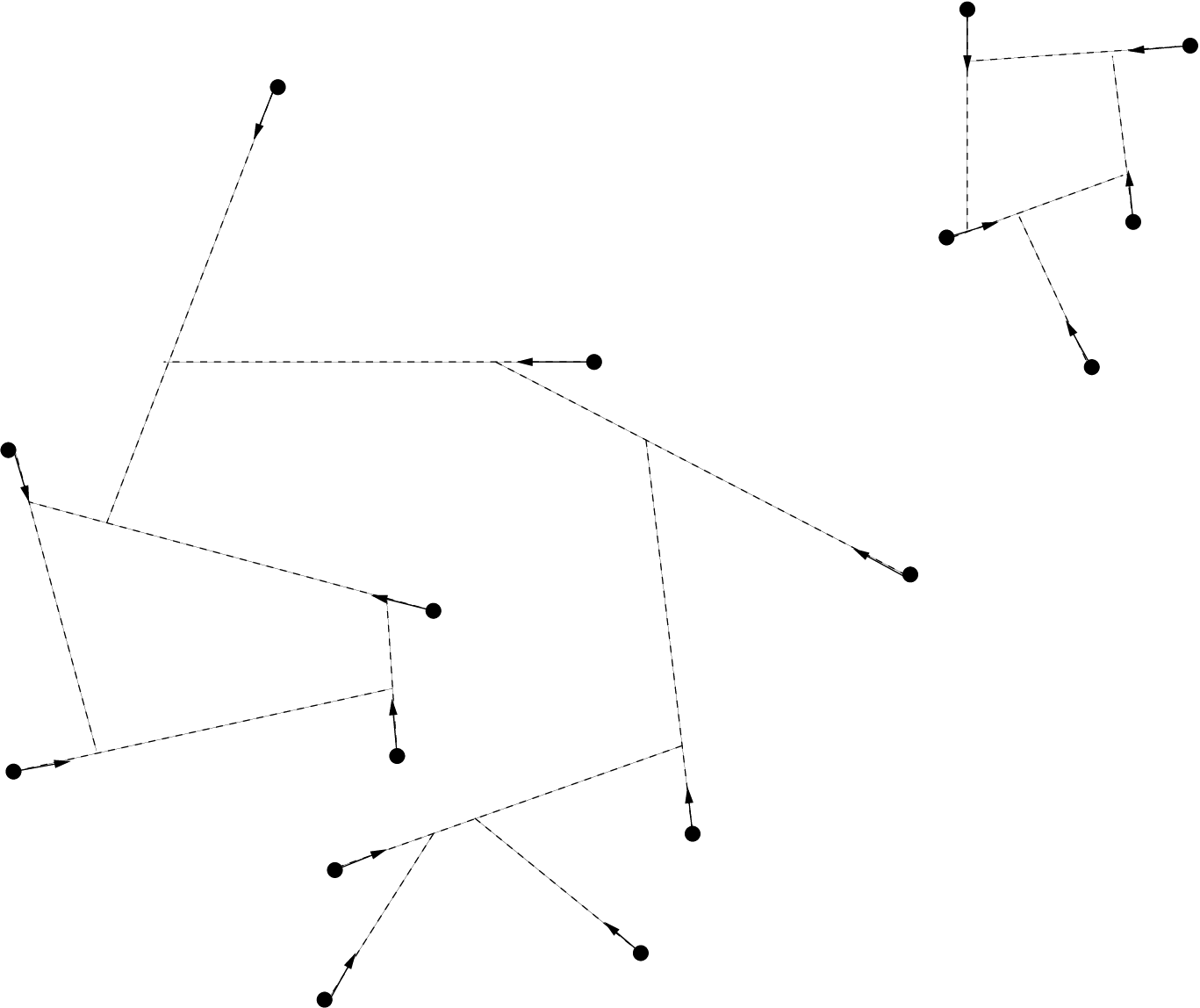}
\hfill
\includegraphics[width=6.5cm]{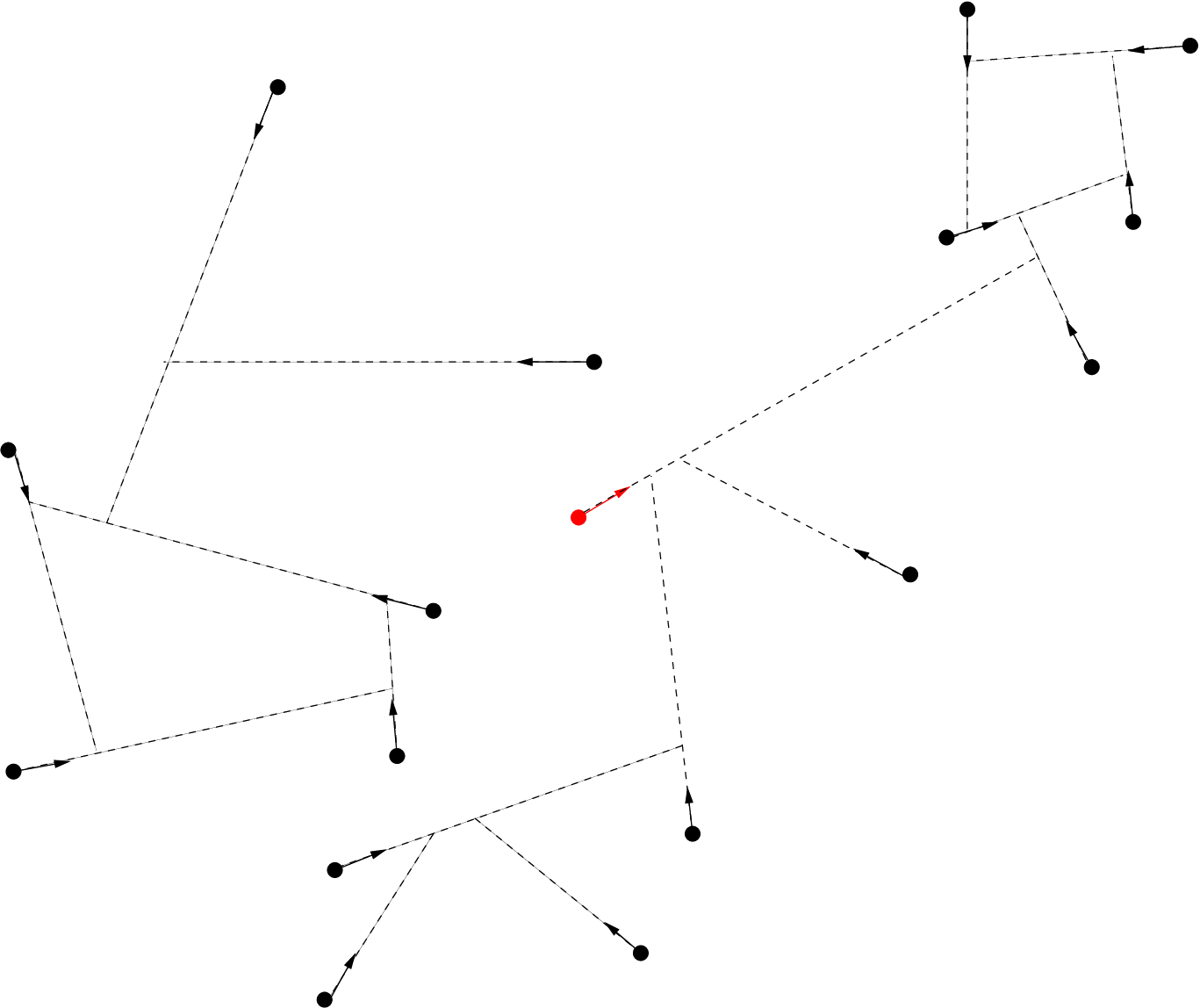}
\caption{In the left part of the picture, we have drawn the geometric graph defined by a finite configuration. Two connected components are obtained. In the right part, we only add one marked point to the configuration. The two initial connected components are modified. Precisely, any point in the initial configuration has a new cluster.}
\label{figure_coup}
\end{center}
\end{figure}

When the growing balls are replaced with growing segments the issue is much more complicated. The two-sided line-segment model is defined via a marked homogeneous Poisson point process $\mathbb{X}$ in $\mathbb{R}^2\times [0;1]$ where the marks are independent and  uniformly distributed on $[0;1]$. At time 0, for any $(\xi,u)\in\mathbb{X}$ a line-segment centred at $\xi$ starts to grow at unit rate in the direction $\pi u$. A one-sided version also exists in which a one-sided segment grows from $\xi$ with direction $2\pi u$. In both models, each line-segment ceases to grow when one of its ends hits another segment. The descending chain argument does not work in this setting. Indeed, when a line-segment hits another one, there is no reason that the hit segment is smaller than the hitting one. Nevertheless the absence of percolation for the two-sided model was conjectured by Daley et al. (Conjecture 7.1 of \cite{daley2014two}) and was proved for the one-sided model in a weaker form (only the four  directions North, East, South and West are allowed) by Hirsch \cite{hirsch2016absence}. Although both references \cite{daley2014two,hirsch2016absence} are rather recent it seems that the natural question of percolation for these line-segments models was known since a while in the probabilistic community. Like any good conjecture, the percolation question for these line-segments models is easy understanding but reveals technical hurdles. On the one hand, any local modification of the marked point process may have huge aftereffects on the final realization of line-segments: see Figure \ref{figure_coup}. On the other hand, the sequence of successively hit line-segments presents no Markovian property or renewal structure.

In this paper, we prove the finite cluster property for the one-sided line-segment model as a consequence of a general result (Theorem \ref{mainTH}, our main result) dealing with outdegree-one graphs. The finite cluster property of the two-sided line segment model can also be obtained following the strategy we have used, the two-sided version does not contain specificities which annihilate the proof.

\begin{figure}[!ht]
\begin{center}
\psfrag{x}{\small{$x$}}
\includegraphics[width=10cm,height=5cm]{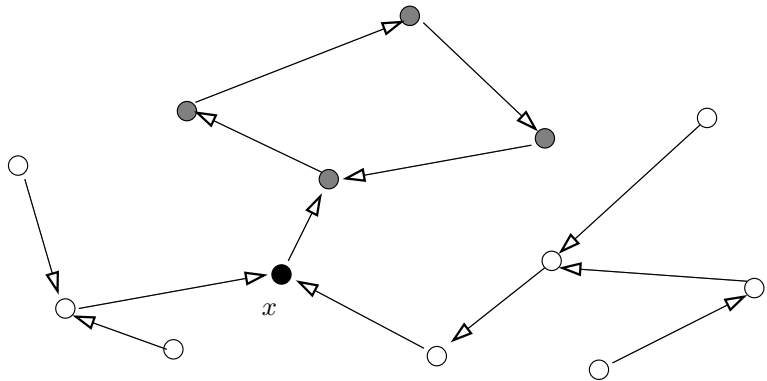}
\caption{\label{fig:cluster} Here is a picture of the cluster of a given vertex $x$. The gray vertices belong to the Forward set of $x$ whereas the white ones are in its Backward set.}
\end{center}
\end{figure}

It is easy (and natural) to interpret the geometric graphs mentioned above as (Poisson) outdegree-one graphs. For the stopped germ-grain models (as lilypond and line-segment models), an oriented edge from $x$ to $y$ is declared when the grain from $x$ hits the grain from $y$. For this reason our main result (Theorem \ref{mainTH}), presented in the general setting of Poisson outdegree-one graph, covers naturally the geometrical setting. Since the graph is oriented, we can define the Forward and Backward sets of any given vertex $x$: see Figure \ref{fig:cluster} for an illustration and Section \ref{sect:Outdegree-oneGraph} for a precise definition. Then, the cluster containing $x$ merely is the union of these both sets.

Thanks to the outdegree-one structure, a forward branch is finite if and only if it contains a loop. By a loop, we mean an finite sequence $x_{1},x_{2},\ldots,x_{n}$ of different vertices for which $x_{i}$ is connected to $x_{i+1}$ for $1\le i\le n-1$ and $x_{n}$ is connected to $x_{1}$. In Figure \ref{fig:cluster}, the Forward set of $x$ contains a loop of size $4$. A general argument for stationary outdegree-one graphs, called mass transport principle, ensures that the absence of forward percolation implies the absence of backward percolation. Henceforth, the aim of our work is to provide general assumptions ensuring that any forward branch of Poisson outdegree-one graphs merges on a loop. This conjecture is supported by the following heuristic (to which our proof strategy does not correspond): there are many loops everywhere in the plane and it is too difficult for an infinite forward branch to avoid all of them.

In \cite{hirsch2016absence}, Hirsch proposes such assumptions in the setting of geometrical graphs. In particular, the author proved a weaker version of the conjecture by Daley et al. for the one-sided line-segment model in which directions are piked uniformly among North, East, South and West. Hirsch's proof consists in stating that any infinite forward branch has to cross a infinite number of ``controlled regions'' in which it can merge on a loop with positive probability and independently of what precedes. To carry out this strategy, Hirsch requires technical assumptions (Section 3 of \cite{hirsch2016absence}) which actually are difficult to check even for the one-sided line-segment model with only four directions. When all directions are allowed, the verification seems likely impossible as says the author himself.

Our proof (of Theorem \ref{mainTH}) differs from the Hirsch's one and deeply exploits the outdegree-one structure of our models. It is based on the following general statement for Poisson outdegree-one graphs which can be roughly interpreted as a counterpart of the mass transport principle: if there exists, with positive probability, an infinite forward branch then the expectation of the size of a typical backward branch is infinite. Our main contribution is to provide minimal assumptions guaranteeing that such expectation is finite, and then ensuring the absence of percolation for a large class of models, containing at least the original line-segment model. As far as we know, this strategy has never been investigated before for proving the absence of percolation in any continuous or discrete models.

Let us describe briefly both assumptions of our main theorem. The first one, called the Loop assumption and denoted by \textbf{(LA)}, assumes that any forward branch merges to a loop if the process is augmented with a finite collection of well-chosen points (without reducing the size of the backward). Roughly speaking, this assumption assures that a loop is possible along a forward branch provided that some points are added. The extra condition on the size of the backward is directly related to the method we use. This condition seems a bit artificial and could be probably relaxed in the future. However we note that it is relatively easy to check in all models we met. The second one, called the Shield assumption and denoted by \textbf{(SA)}, is directly inspired from the ones by Hirsch. More or less, it assumes that with high probability, the graph contains no edge crossing large boxes. 

{Moreover, the conclusion of Theorem \ref{mainTH}, i.e. the absence of percolation, does not hold if only one hypothesis among \textbf{(LA)} and \textbf{(SA)} is satisfied. Whereas the Loop assumption appears as an essential property to prevent percolation, it is not clear whether the Shield assumption is really needed. We construct an outdegree-one graph which satisfies only the Loop assumption and percolates.}

As mentioned before our main application is the absence of percolation for the line-segment model introduced by Daley et al. We investigate also another model which is inspired by the geometrical navigation defined in \cite{bonichon2011asymptotics}. See Theorem \ref{theo:models}.

The paper is organized as follows. In Section 2, we provide a precise description of Poisson outdegree-one models and give examples. In Section 3, we formulate our two assumptions and the main result (Theorem \ref{mainTH}) ensuring the absence of percolation. Section 4 is devoted to its proof and finally, in Section 5, we check that both models introduced in Section 2 satisfy the assumptions of Theorem \ref{mainTH}.

\section{General model and examples}

\subsection{Notations}

All the models of this paper take place in the Euclidean space $\Rd$. The configuration space $\C$ on $\Rd$ with marks in $[0;1]$ is defined by
$$
\C = \Big\{ \varphi\subset\Rd\times[0;1] \,;\, N_{\Lambda}(\varphi)<\infty, \text{ for any bounded } \Lambda\subset\Rd \Big\}
$$
where $N_{\Lambda}(\varphi)=\#(\varphi\cap(\Lambda\times[0;1]))$ denotes the number of marked points of $\varphi$ whose first ordinate lies in $\Lambda$. Any other choice of compact set for the marks could be considered with slight modifications in the following. Let us denote by $\varphi_{germs}$ the projection of any given configuration $\varphi\in\mathscr{C}$ onto $\Rd$: $\varphi_{germs}=\{\xi\,;\,(\xi,\cdot)\in\varphi\}$. For a given subset $\Lambda$ of $\Rd$, and $\varphi\in\mathscr{C}$, $\varphi_{\Lambda}$ denotes the set of points of $\varphi$ included in $\Lambda\times [0;1]$: $\varphi_{\Lambda}=\lbrace(\xi,u)\in\varphi\ ;\ \xi\in\Lambda\rbrace$.

As usual, the configuration space $\mathscr{C}$ is equipped with the $\sigma$-algebra
$$
\mathscr{S} = \sigma \Big( P_{(A,n)} \,;\, A \,\mbox{ Borel set of }\, \Rd\times[0;1]\text{ a Borel set} ,\, n\ge 0 \Big) ~,
$$
generated by the counting events $P_{(A,n)}=\{\varphi\in\mathscr{C} ; \#(\varphi\cap A)\le n \}$. Similarly, for a any subset $\Lambda\subset\Rd$, we define the $\sigma$-algebra of events in $\Lambda$ by 
$$
\mathscr{S}_{\Lambda} = \sigma \Big( P_{(A,n)} \,;\, A \,\mbox{ Borel set of }\, \Lambda\times[0;1] ,\, n\ge 0 \Big) ~.
$$
Let $v\in\Rd$. The translation operator $\tau_v$ acts on $\Rd$, $\Rd\times[0;1]$ and $\mathscr{C}$ as follows: for any $w\in\Rd$, $x=(\xi,u)\in\Rd\times[0;1]$ and $\varphi\in\C$, we set $\tau_v(w)=v+w$, $\tau_v(x)=(\xi+v,u)$ and $\tau_v(\varphi)=\cup_{x\in\varphi}\lbrace \tau_{v}(x)\rbrace$. Finally, a subset $\C'\subset\C$ is said \textit{translation invariant} whenever $\tau_v(\C')\subset\C'$, for any vector $v\in\Rd$.

\subsection{The outdegree-one model}
\label{sect:Outdegree-oneGraph}

In our setting, an \textit{outdegree-one graph} is an oriented graph whose vertex set is given by a configuration $\varphi\in\C$ and having exactly one outgoing edge per vertex. Such graph can be described by a \textit{graph function} which determines, for any vertex, its outgoing neighbour. Note that the marks in $[0;1]$ will be used as random contributions to this connection mechanism. See examples in Section \ref{sect:examples}.

\begin{defin}
\label{DefinitionGM}
Let $\C'\subset\C$ be a translation invariant set. A function $h$ from $\C'\times(\Rd\times[0;1])$ to $\Rd\times[0;1]$ is called a \textbf{graph function} if:
\begin{itemize}
\item[(i)] $\forall\varphi\in\C', \forall x\in\varphi , \, h(\varphi,x)\in\varphi\backslash\{x\}$;
\item[(ii)] $\forall v\in\Rd, \forall\varphi\in\C', \forall x\in\varphi, \, h(\tau_{v}(\varphi),\tau_{v}(x)) = \tau_{v}(h(\varphi,x))$.
\end{itemize}
The couple $(\C',h)$ is then called an \textbf{outdegree-one model}.
\end{defin}

In the sequel, let us consider an outdegree-one model $(\C',h)$ and a configuration $\varphi\in\C'$. The associated graph is made up of edges $(x, h(\varphi,x))$, for all $x\in\varphi$. In dimension $d=2$, such graphs are not necessarily planar.

Let us describe the clusters of this graph. Let $x\in\varphi$. The \textit{Forward set} $\text{For}(x,\varphi)$ of $x$ in $\varphi$ is defined as the sequence of the outgoing neighbours starting at $x$:
$$
\text{For}(x,\varphi) = \{ x, h(\varphi,x), h(\varphi, h(\varphi,x)), \ldots \} ~.
$$
The Forward set $\text{For}(x,\varphi)$ is a branch of the graph, possibly infinite. The \textit{Backward set} $\text{Back}(x,\varphi)$ of $x$ in $\varphi$ contains all the vertices $y\in\varphi$ having $x$ in their Forward set:
$$
\text{Back}(x,\varphi) = \{ y\in\varphi \,;\, x \in \text{For}(y,\varphi) \} ~.
$$
The Backward set $\text{Back}(x,\varphi)$ admits a tree structure whose $x$ is the root. The Forward and Backward sets of $x$ may overlap; they (at least) contain $x$. Their union forms the \textit{Cluster} of $x$ in $\varphi$:
$$
C(x,\varphi) = \text{For}(x,\varphi) \cup \text{Back}(x,\varphi) ~.
$$
The \textit{Cluster} $C(x,\varphi)$ is a subset of the connected component of $x$ in $\varphi$, the absence of infinite cluster in a given outdegree-one graph is nothing else than the absence of infinite connected component.

Our main theorem (Theorem \ref{mainTH}) states that for a large class of random models, all the clusters are a.s. finite. In particular, it is not difficult to observe that the Forward set $\text{For}(x,\varphi)$ is finite if and only if it contains a \textit{loop}, i.e. a subset $\{y_{0},\ldots,y_{l-1}\}\subset\text{For}(x,\varphi)$, with $l\ge2$, such that for any $0\le i\le l-1$, $h(\varphi,y_{i})=y_{i+1}$ (where the index $i+1$ is taken modulo $l$). In this case, the integer $l$ is called the \textit{size} of the loop. Furthermore, the outdegree-one property implies that there is at most one loop in a cluster. Hence, a finite cluster is made up of one loop with some finite trees rooted at vertices of the loop (see Figure \ref{fig:cluster}). Obviously, this notion of loops will be central in our study.

\subsection{Random outdegree-one model}

Let $Q$ be a probability measure on $[0;1]$ such $Q(\theta)>0$ for any open set $\theta$ in $[0;1]$ and let us denote by $\lambda_{d}$ the Lebesgue measure on $\Rd$. We consider a Poisson point process (PPP) $\mathbb{X}$ on $\C$ with intensity $\lambda_{d}\otimes Q$. This means that the random variable $\#(\mathbb{X}\cap A)$ follows a Poisson distribution with parameter $\lambda_{d}\otimes Q(A)$, for any bounded Borel set $A\subset\Rd\times [0;1]$. By a standard change of scale, any other (stationary) intensity measure of the form $z\lambda_{d}\otimes Q$ with $z>0$ could be considered. The process $\mathbb{X}$ can also be interpreted as a stationary PPP on $\Rd$ with intensity one in which all the Poisson points are independently marked with distribution $Q$, and such that the marks are also independent of the locations of the Poisson points. For the two examples considered in this work (see Section \ref{sect:examples}), $Q$ is the uniform distribution on $[0;1]$. But other mark distributions could be foreseen: this will be discussed in Section \ref{sect:TH}.

Finally, let us denote by $(\Omega,\mathscr{F},\mathbb{P})$ a probability space on which the PPP  $\mathbb{X}$ is defined.

\begin{defin}
\label{defin:RGM}
Let $(\C',h)$ be an outdegree-one model. If $\mathbb{P}(\mathbb{X}\in\C')=1$ then the triplet $(\C',h,\mathbb{X})$ is called a \textbf{random outdegree-one model}.
\end{defin}

\subsection{Two examples}
\label{sect:examples}

Let $Q$ be the uniform distribution on $[0;1]$.

\subsubsection{The line-segment model}

Our first model is a unilateral version of the model studied in \cite{daley2014two} (called \textit{Model 1} therein) and mentioned in Section \ref{SectionIntro}. It is also a generalization of the model studied in \cite{hirsch2016absence}.

The \textit{line-segment model} is based on a stopping germ-grain protocol defined as follows. Let us consider a marked configuration $\varphi$ in $\Rde\times[0;1]$. At the same time (say $t=0$), for any marked point $(\xi,u)\in\varphi$, an half line-segment starts growing (at unit rate) from $\xi$ according to the direction $2\pi u$. Each line-segment ceases to grow whenever its end point hits another line-segment. But the stopping one continues its growth.

Let us denote by $\C'$ the configuration space for which the above dynamic is well defined. This means that each line-segment is eventually stopped by exactly one other line-segment. Of course, the set $\C'$ is translation invariant. We can then define a graph function $h$ encoding the line-segment model: given $\varphi\in\C'$ and $x\in\varphi$, the image $h(\varphi,x)$ refers to the stopping line-segment of $x$. This construction clearly provides an outdegree-one graph.

The authors in \cite{daley2014two} proved that $\mathbb{P}(\mathbb{X}\in\C')=1$. Roughly speaking, they proved that for almost all configuration, the unique stopping segment of any point can be determined by a finite algorithm. Hence they have checked the existence of the two-sided and one-sided line segment model. Therefore, according to Definition \ref{defin:RGM}, $(\C',h,\mathbb{X})$ is a random outdegree-one model.

\subsubsection{The navigation model}

\begin{figure}[!ht]
\begin{center}
\includegraphics[width=8.5cm,height=5cm]{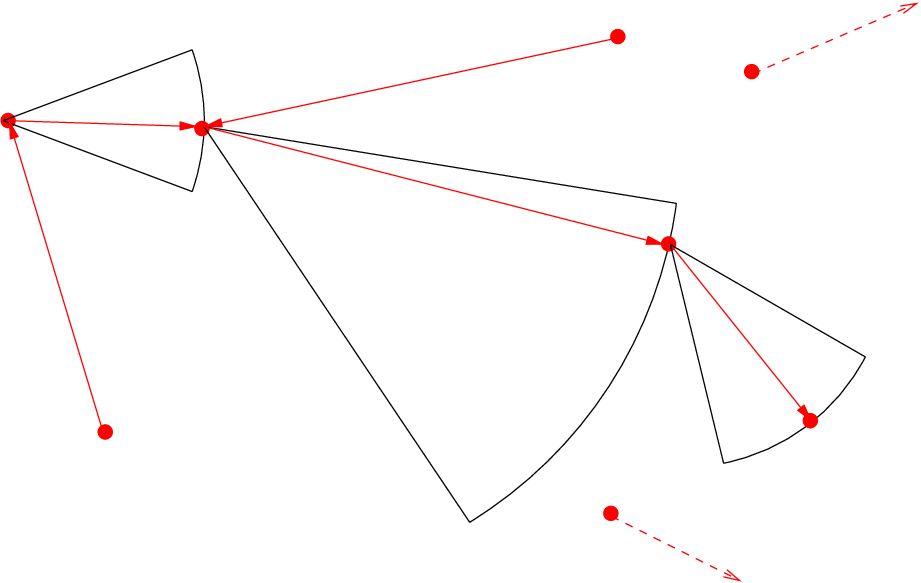}
\caption{\label{fig:Navigation} Here is a scheme of the navigation model. Remark that in the navigation model, two marked points suffice to form a loop.}
\end{center}
\end{figure}

Let us define the \textit{navigation model} introduced in \cite{bonichon2011asymptotics}. Let $\epsilon\in(0;\pi]$ be an extra parameter and $\varphi$ be a configuration in $\Rde\times[0;1]$. Each marked point $x=(\xi,u)\in\varphi$ defines a semi-infinite cone with apex $\xi$, direction $2\pi u$ and opening angle $\epsilon$:
$$
C(x) = \xi + \Big\{ (r\cos(\alpha),r\sin(\alpha)) \,;\, r>0 \,\mbox{ and }\, |\alpha-2\pi u|<\epsilon \Big\} ~.
$$
In the navigation model, each marked point $x=(\xi,u)$ is connected to $(\xi',\cdot)\in\varphi$ where $\xi'$ is the closest element to $\xi$ living in the cone $C(x)$, i.e. such that
\begin{equation}
\label{defnavigation}
\|\xi-\xi'\|_{2} = \min \Big\{ \|\xi-\xi''\|_{2} \,;\, \xi'' \in C(x) \cap \varphi_{germs} \Big\} ~.
\end{equation}
See Figure \ref{fig:Navigation}. Of course this connection procedure produces an outdegree-one graph provided those closest elements exist and are unique. Let us denote by $\C'$ the set of such configurations. This is a translation invariant set. In this setting, $h(\varphi,x)$ is defined as the unique marked point $(\xi',\cdot)$ where $\xi'$ satisfies (\ref{defnavigation}).

Using standard properties of the PPP, it is easy to show that $\mathbb{P}(\mathbb{X}\in\C')=1$ and therefore $(\C',h,\mathbb{X})$ is a random outdegree-one model.

In the case where $\epsilon=\pi$, the marks have no longer importance and the navigation model actually coincides with the nearest neighbour graph.

\section{Results}
\label{sect:TH}

We first establish in Theorem \ref{mainTH} the absence of percolation for all random outdegree-ones $(\C',h,\mathbb{X})$ satisfying two general assumptions, namely the \textit{Loop and Shield assumptions}, which are described and commented below. Thus, Theorem \ref{theo:models} asserts that the line-segment model and the navigation model verify these two assumptions and therefore do not percolate.\\

\noindent
{\textbf{Loop assumption}}

The Loop assumption \textbf{(LA)} mainly expresses the possibility for any marked point $x\in\varphi$ to break its Forward set by adding a finite sequence of marked points $(x_1,\ldots,x_k)$ to the current configuration $\varphi$.

Let $\varphi\in\C'$ and $k$ be a positive integer. The configuration $\varphi$ is said $k$\textit{-looping} if for any $x\in\varphi$, there exists an open ball $A_{x}\subset(\Rd\times[0;1])^{k}$ such that, for all $(x_{1},\ldots,x_{k})\in A_{x}$:
\begin{itemize}
\item[(i)] $\text{For}(x,\varphi\cup\{x_{1},\ldots,x_{k}\}) \subset \{x,x_{1},\ldots,x_{k}\}$,
\item[(ii)] $\forall 1\le i\le k $, $x_{i}$ belongs to the connected component of $x$ in $\varphi\cup\{x_{1},\dots,x_{k}\}$ and
$$\text{For}(x_{i},\varphi\cup\{x_{1},\dots,x_{k}\})\subset\{x,x_{1},\dots,x_{k}\}, $$
\item[(iii)] {There exists a positive constant $C_1$ which does not depend on $x$ such that
$$
\# \text{Back}(x,\varphi\cup\{x_{1},\ldots,x_{k}\}) \geq \# \text{Back}(x,\varphi) - C_{1}
$$
(in $\mathbb{N}\cup\{\infty\}$).}
\end{itemize}


Given $x$, the three conditions above can be interpreted as a local modification of the configuration $\varphi$ which breaks the Forward set of $x$ without decreasing {too much} the cardinality of its Backward set. Item $(i)$ is very natural to obtain a finite cluster. Item $(ii)$-- combined with \textbf{(SA)} --will be crucial {to localize the added points $x_i$'s around $x$ and then} to guarantee the construction of local events in Section \ref{sect:InfiniteAlmostLooping}. Item $(iii)$ is more technical and will appear in the proof of Proposition \ref{prop:FromAlmostLoopingToLooping}. {Let us note that for the line-segment model and the navigation model, we will prove that $\text{Back}(x,\varphi)$ is included in $\text{Back}(x,\varphi\cup\{x_{1},\ldots,x_{k}\})$ which implies $(iii)$ with $C_1=0$.} Finally, the choice of the integer $k$ will be adapted to the random outdegree-one model $(\C',h,\mathbb{X})$.

We will say that the random outdegree-one model $(\C',h,\mathbb{X})$ satisfies \textbf{(LA)} if there exists a positive integer $k$ such that
$$
\mathbb{P}(\mathbb{X}\, \mbox{  is $k$-looping } ) = 1 ~.
$$

\noindent
{\textbf{Shield assumption}}

The Shield assumption \textbf{(SA)} is a kind of strong stabilizing property for the random outdegree-one $(\C',h,\mathbb{X})$ and has been first introduced in a slightly different way in \cite{hirsch2016absence}.

\begin{figure}[!ht]
\begin{center}
\psfrag{a}{\small{$m \mathcal{A}_{1}$}}
\psfrag{b}{\small{$m V$}}
\psfrag{c}{\small{$m \mathcal{A}_{2}$}}
\psfrag{d}{\small{$\xi$}}
\includegraphics[width=8cm,height=6.5cm]{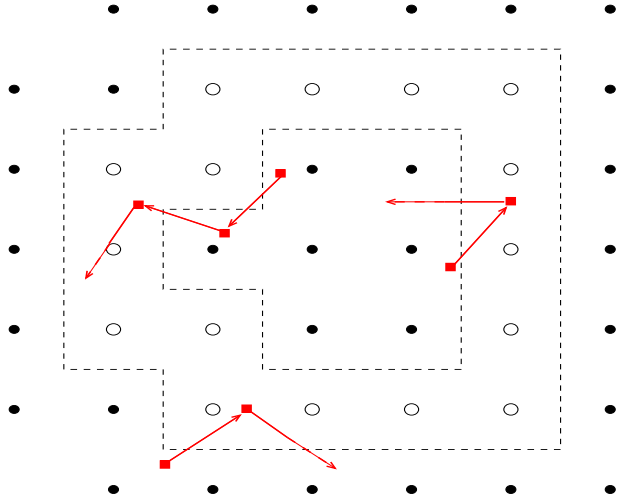}
\caption{\label{fig:SA} The white points are the elements of $mV$ while the black ones are those of $m\mathcal{A}_{1}$ (inside $mV$) and $m\mathcal{A}_{2}$ (outside $mV$). Red squares are points of $\varphi$. If $\tau_{-mz}(\varphi)\in\mathscr{E}_{m}$ for all $z\in V$, then it is impossible for a (red) segment $[\xi;\xi']$, where $x=(\xi,\cdot)$ and $h(\varphi,x)=(\xi',\cdot)$, to cross the set $m(V\oplus[-\alpha;\alpha]^{d})$ from $m\mathcal{A}_{1}$ to $m\mathcal{A}_{2}$-- or from $m\mathcal{A}_{2}$ to $m\mathcal{A}_{1}$ by symmetry of Equation (\ref{conditionSA}) w.r.t. indexes $1$ and $2$.}
\end{center}
\end{figure}

We will say that the random outdegree-one $(\C',h,\mathbb{X})$ satisfies \textbf{(SA)} if there exist a positive integer $\alpha$ and a sequence of events $(\mathscr{E}_{m})_{m\ge 1}$ such that:
\begin{itemize}
\item[(i)] $\forall m\ge 1$, $\mathscr{E}_{m}\in\mathscr{S}_{[-\alpha m;\alpha m]^{d}}$;
\item[(ii)] $\mathbb{P}(\mathscr{E}_{m}) \underset{m\to\infty}{\longrightarrow} 1$;
\item[(iii)] Consider the lattice $\Z^{d}$ with edges given by $\{\{z,z'\}, |z-z'|_{\infty}=1\}$ and any three disjoint subsets $A_{1}, A_{2}, V$ of $\Z^{d}$ such that $\forall i=1,2$, the boundary $\partial A_{i}=\{z\in\Z^{d}\backslash A_{i}, \exists z'\in A_{i}, |z-z'|_{\infty}=1\}$ is included in $V$. Let us set
$$
\mathcal{A}_{i} = \Big( A_{i}\oplus \Big[ -\frac{1}{2};\frac{1}{2} \Big]^{d} \Big) \setminus (V\oplus[-\alpha;\alpha]^{d}) ~.
$$
Then, for all $m$ and for any configurations $\varphi,\varphi'\in\C'$ such that $\tau_{-mz}(\varphi)\in\mathscr{E}_{m}$ for all $z\in V$, the following holds:
\begin{equation}
\label{conditionSA}
\forall x \in \varphi_{m\mathcal{A}_{1}} , \,  h(\varphi,x) = h(\varphi_{m\mathcal{A}_{2}^{c}} \cup \varphi'_{m\mathcal{A}_{2}},x) ~.
\end{equation}
\end{itemize}
In Condition $(iii)$, the set $mV$ acts as an uncrossable obstacle, i.e. a shield, between sets $m\mathcal{A}_{1}$ and $m\mathcal{A}_{2}$. See Figure \ref{fig:SA}. Equation (\ref{conditionSA}) says that the outgoing neighbour of any $x\in\varphi_{m\mathcal{A}_{1}}$ does not depend on the configuration on $m\mathcal{A}_{2}$. In particular, $h(\varphi,x)\in \varphi_{m\mathcal{A}_{2}^{c}}$.

Here is our main result.

\begin{theo}
\label{mainTH}
Any random outdegree-one $(\C',h,\mathbb{X})$ satisfying \textbf{(LA)} and \textbf{(SA)} does not percolate with probability $1$:
$$
\mathbb{P}(\forall x \in \mathbb{X} , \, \#C(x,\mathbb{X}) < \infty ) = 1 ~.
$$
\end{theo}

Theorem \ref{mainTH} is proved in Section \ref{SectionPMT} and a sketch is given in Section \ref{sect:sketch}.\\
Checking that both models of Section \ref{sect:examples} satisfy \textbf{(LA)} and \textbf{(SA)} (this is done in Section \ref{SectionPC}), we get:

\begin{theo}
\label{theo:models}
The line-segment model and the navigation model do not percolate with probability $1$.
\end{theo}

The choice of a suitable integer $k$ such that the random outdegree-one $(\C',h,\mathbb{X})$ satisfies \textbf{(LA)} actually depends on the model. For example, $3$ marked points suffice to make a loop for the line-segment model and only $1$ for the navigation model: see respectively Propositions \ref{prop:CyclingAssGrowing} and \ref{Prop:CycAssNav}.

Any variant of the line-segment model or the navigation model in which the uniform distribution $Q$ is replaced with a probability measure $Q'$ defined by
$$
Q'(A) = \int_{A} f(x) Q(dx) ~,
$$
where $f$ is a positive function on $[0;1]$, has to satisfy \textbf{(LA)} and \textbf{(SA)}.

Furthermore, the law of the size of a typical cluster for a given model satisfying both assumptions \textbf{(LA)} and \textbf{(SA)} would be an interesting further result. It is also possible to investigate asymptotics of the number of different clusters on expanding windows. In the case of the lilypond model some results of this kind have been established in \cite{last2013percolation}.\\

{Let us end this section with discussing the need for assumptions \textbf{(LA)} and \textbf{(SA)} in Theorem \ref{mainTH}. To do it, let us consider the Directed Spanning Forest (DSF) in $\mathbb{R}^2$ with direction $(1,0)$ in which each Poisson point $x$ in $\mathbb{X}$ is connected to its nearest Poisson point $h(\mathbb{X},x)$ having a larger abscissa, i.e.
$$
h(\mathbb{X},x) := \mbox{argmin} \big\{ \|x-y\|_{2} : y \in \mathbb{X} \,\mbox{ and }\, \langle y-x , (1,0) \rangle \geq 0 \big\} ~,
$$
where $\langle \cdot , \cdot \rangle$ denotes the inner product. See \cite{CT} for the study of infinite branches of the DSF.}

{It is not difficult to check that the DSF is a random outdegree-one model-- whithout marks --satisfying \textbf{(SA)} but not \textbf{(LA)}, and percolating since by construction a semi-infinite branch starts at each vertex.}

{Now, we are going to modify the DSF into a new random outdegree-one model, say DSF$^{\ast}$ corresponding to a graph function $h^{\ast}$, satisfying this time \textbf{(LA)} but not \textbf{(SA)}, and still percolating. Let $w:\mathbb{N}\to(0,+\infty)$ be a decreasing function. For any couple of vertices $(x,h(\mathbb{X},x))$ such that
\begin{itemize}
\item[(a)] $x$ and $h(\mathbb{X},x)$ are mutually (Euclidean) nearest neighbours;
\item[(b)] $\{y\in\mathbb{X}: h(\mathbb{X},y)=h(\mathbb{X},x)\}=\{x\}$;
\item[(c)] $\|x-h(\mathbb{X},x)\|_{2} \leq w(\# \text{Back}(x,\mathbb{X}))$;
\end{itemize}
then $h^{\ast}(\mathbb{X},x)=h(\mathbb{X},x)$ and $h^{\ast}(\mathbb{X},h(\mathbb{X},x))=x$ (Recall that in the DSF the Backward set $\text{Back}(x,\mathbb{X})$ of any vertex $x$ is a.s. finite : \cite{CT}, Theorem 8). Otherwise, outgoing edges remain unchanged. In other words, the only change between the DSF and DSF$^{\ast}$ consists in the possibility to create a loop of size $2$ with the mutual nearest neighbours $x$ and $h(\mathbb{X},x)$ satisfying Items (b) and (c). First remark that the resulting graph DSF$^{\ast}$ satisfies \textbf{(LA)} with $k=1$. Indeed, it is always possible to add a point $x_1$ as close as we want to $x$ such that:
\begin{enumerate}
\item $\langle x_1 - x , (1,0) \rangle \geq 0$.
\item $x$ and $x_1$ are mutually nearest neighbours in $\mathbb{X}\cup\{x_1\}$.
\item $\| x - x_1 \|_{2} \leq w(\# \text{Back}(x,\mathbb{X}))$.
\item The add of $x_1$ does not reduce the Backward set of $x$: $\#\text{Back}^{\ast}(x,\mathbb{X}) \leq \#\text{Back}^{\ast}(x,\mathbb{X}\cup\{x_1\})$.
\end{enumerate}
Items 1 and 2 ensure that $h(\mathbb{X}\cup\{x_1\},x)=x_1$. Combined with Item 3, we get that $(x,x_1)$ forms a loop in DSF$^{\ast}$. So Items $(i)$ and $(ii)$ of \textbf{(LA)} are clearly true. In the case where $x$ and $h^{\ast}(\mathbb{X},x)$ form a loop in DSF$^{\ast}$, the add of $x_1$ makes $h^{\ast}(\mathbb{X},x)$ (and only it by Item (c)) come out of the Backward set of $x$. But this output is offset by the input of $x_1$, which leads to Item 4 and then to Item $(iii)$ of \textbf{(LA)}. Note also that \textbf{(SA)} does not hold for the DSF$^{\ast}$ since we possibly need to explore all the set $\text{Back}(x,\mathbb{X})$ to determine if $x$ and $h(\mathbb{X},x)$ form a loop in DSF$^{\ast}$. It then remains to show that DSF$^{\ast}$ admits infinite paths. Let $(x_n)_{n\geq0}$ be the forward path (in the DSF) starting at a given vertex $x=x_0\in\mathbb{X}$, with $h(\mathbb{X},x_n)=x_{n+1}$ for any $n$. Let us choose the decreasing function $w$ such that, for any $n$, $\mathbb{P}(\|x_n-x_{n+1}\|_{2} \leq w(n))<2^{-(n+1)}$. Then, with positive probability, $x$ admits also an infinite forward path for DSF$^{\ast}$. Indeed,
\begin{eqnarray*}
\mathbb{P} \big( \# \text{For}^{\ast}(x,\mathbb{X}) < \infty \big) & \leq & \mathbb{P} \big( \exists n\geq 0 , \| x_{n}-x_{n+1} \|_{2} \leq w(\# \text{Back}(x_{n},\mathbb{X})) \big) \\
& \leq & \sum_{n\geq 0} \mathbb{P} \big( \| x_{n}-x_{n+1} \|_{2} \leq w(n) \big) \; < \; 1 ~.
\end{eqnarray*}
Hence, among all the (necessarily infinite) forward branches of the DSF, some of them remain unchanged (and then infinite too) when passing from DSF to DSF$^{\ast}$.}

\section{Proof of Theorem \ref{mainTH}}
\label{SectionPMT}

\subsection{Sketch of the proof}
\label{sect:sketch}

We have to prove that any random outdegree-one model satisfying \textbf{(LA)} and \textbf{(SA)} does not contain any infinite cluster with probability $1$, i.e.
\begin{equation}
\label{goalMainTH}
\mathbb{P} \big( \forall x\in\mathbb{X} , \, \#\text{For}(x,\mathbb{X})<\infty \, \mbox{ and } \, \#\text{Back}(x,\mathbb{X})<\infty \big) = 1 ~.
\end{equation}
First, thanks to a standard mass transport argument (Proposition \ref{PropTM} of Section \ref{sect:ForwardtoBackward}), we can reduce the proof of the absence of percolation to the absence of forward percolation:
\begin{equation}
\label{ForwardtoBackward}
\mathbb{P} \big( \forall x\in\mathbb{X} , \, \#\text{For}(x,\mathbb{X}) < \infty \big) = 1 \; \Longrightarrow \; \mathbb{P} \big( \forall x\in\mathbb{X} , \, \#\text{Back}(x,\mathbb{X}) < \infty \big) = 1 ~.
\end{equation}
Then, we proceed by contradiction and assume that, with positive probability, an infinite forward branch starts at a typical marked point $\Theta$:
\begin{equation}
\label{InfiniteForwardTypical}
\mathbb{P} \big( \#\text{For}(\Theta,\mathbb{X}_{\Theta}) = \infty \big) > 0 ~,
\end{equation}
where $\mathbb{X}_{\Theta}$ denotes the configuration $\mathbb{X}\cup\{\Theta\}$ with $\Theta=(0,U)$ and $U$ is an uniform random variable in $[0;1]$. Two central notions here are \textit{looping points} and \textit{almost looping points}. To sum up, a looping point admits a finite Forward set, i.e. a forward branch ending with a loop. An almost looping point is set to become a looping point by adding some suitable marked points. See respectively Definitions \ref{defi:LoopingPoint} and \ref{DefinitionALP}. In Section \ref{sect:InfiniteAlmostLooping}, we use intensively hypotheses \textbf{(LA)} and \textbf{(SA)} to prove that (\ref{InfiniteForwardTypical}) forces the infinite branch $\text{For}(\Theta,\mathbb{X}_{\Theta})$ to contain an infinite number of almost looping points:
\begin{equation}
\label{InfiniteAlmostLooping}
\mathbb{P} \big( \# \{ y \in \text{For}(\Theta,\mathbb{X}_{\Theta}) ;\, y \text{ is an almost looping point of } \, \mathbb{X}_{\Theta} \} = \infty \big) > 0 ~.
\end{equation}
This is Proposition \ref{PropLiggett}. Heuristically, such event should not occur since it produces an infinite number of opportunities to break the forward branch by adding points. Precisely, Proposition \ref{prop:contradiction} allows to convert the forward result (\ref{InfiniteAlmostLooping}) to a backward one:
\begin{equation}
\label{ExpectInfiniteAlmost}
\mathbb{E} \big[ \#\text{Back}(\Theta,\mathbb{X}_{\Theta}) \1_{\{\Theta \text{ is an almost looping point of } \mathbb{X}_{\Theta} \}} \big] = \infty ~. 
\end{equation}
Thus, by adding some suitable marked points, (\ref{ExpectInfiniteAlmost}) implies that the mean size of the Backward set of a typical looping point is infinite (Proposition \ref{prop:FromAlmostLoopingToLooping}):
\begin{equation}
\label{ExpectInfiniteLooping}
\mathbb{E} \big[ \#\text{Back}(\Theta,\mathbb{X}_{\Theta}) \1_{\{\Theta \text{ is a looping point of } \mathbb{X}_{\Theta} \}} \big] = \infty ~. 
\end{equation}
This actually is the only place where the condition $(iii)$ of \textbf{(LA)} is used. Finally, another use of the mass transport principle (Proposition \ref{PropBLP}) makes statement (\ref{ExpectInfiniteLooping}) impossible. This contradiction achieves the proof of Theorem \ref{mainTH}.


\subsection{Absence of backward percolation}
\label{sect:ForwardtoBackward}

Using the mass transport principle (Lemma \ref{lemmetransport}), we show that the backward percolation is impossible whenever the forward percolation does not occur. This standard argument is formulated in \cite{benjamini1999critical} and used p.18 of \cite{last2013class} and p.4 of \cite{hirsch2016absence}. 

\begin{prop}
\label{PropTM}
The following implication holds:
$$
\mathbb{P} \big[ \forall x\in\mathbb{X} , \, \#\text{For}(x,\mathbb{X})<\infty \big] = 1 \; \Longrightarrow \; \mathbb{P} [ \forall x\in\mathbb{X} , \, \#\text{Back}(x,\mathbb{X})<\infty ] = 1 ~.
$$
\end{prop}

Let us consider a configuration $\varphi\in\C'$ and a bounded subset $\Lambda$ of $\Rd$. A marked point $x\in\varphi$ is said \textit{looping inside} $\Lambda$ (for $\varphi$) if its Forward set $\text{For}(x,\varphi)$ contains a loop $\{y_{1},\ldots,y_{l}\}$,  and the center of mass of the set $\{y_{1},\ldots,y_{l}\}$ belongs to $\Lambda$ (we will also say that the entire forward set is looping inside $\Lambda$). Thus, for any $z\in\Z^{d}$, let us define the following set:
$$
Q_{z}(\varphi) = \Big\{ x \in \varphi ;\, x \text{ is looping inside } z \oplus\Big[-\frac{1}{2};\frac{1}{2}\Big]^{d} \Big\} ~.
$$

\begin{lem}
\label{lemmetransport}
Let $z\in\Z^{d}$. Then, $\mathbb{E}[\#Q_{z}(\mathbb{X})]<\infty$.
\end{lem}

The proof of Lemma \ref{lemmetransport} is based on the mass transport principle.

\begin{dem}[of Lemma \ref{lemmetransport}.]
By stationarity, it is enough to prove that $\mathbb{E}[\#Q_{0}(\mathbb{X})]$ is finite.  For $y,z\in\Z^{d}$, we denote by $Q^{y}_{z}(\mathbb{X})$ the elements $x\in Q_{z}(\mathbb{X})$ whose first ordinate is in $y\oplus\Big[-\frac{1}{2};\frac{1}{2}\Big]^{d}$. Then,
\begin{eqnarray*}
\mathbb{E}[\#Q_{0}(\mathbb{X})] & = & \sum_{y\in\Z^{d}} \mathbb{E}[ \#Q^{y}_{0}(\mathbb{X})] \\
& = & \sum_{y\in\Z^{d}} \mathbb{E}[ \#Q^{0}_{-y}(\mathbb{X})]
\end{eqnarray*}
where the latter equality is due to the stationarity of the PPP $\mathbb{X}$ and the graph function $h$. Now, each cluster in $\mathbb{X}$ a.s. contains at most one loop. This means that
$$
\text{a.s. } \, \sum_{y\in\Z^{d}} \#Q^{0}_{-y}(\mathbb{X}) \leq N_{[-\frac{1}{2},\frac{1}{2}]^{d}}(\mathbb{X}) ~.
$$ 
Since the PPP $\mathbb{X}$ has intensity $1$, it follows $\mathbb{E}[\#Q_{0}(\mathbb{X})]\leq 1$.
\end{dem}

The proof of Proposition \ref{PropTM} is an immediate consequence of Lemma \ref{lemmetransport}.

\begin{dem}[of Proposition \ref{PropTM}.]
Let us assume that with positive probability there exists $x\in\mathbb{X}$ whose Backward set is infinite. By hypothesis, its Forward set a.s. contains a loop. So, we can find a deterministic $z\in\Z^{2}$ such that
$$
\mathbb{P} \Big[ \exists x \in \mathbb{X} ;\, \#\text{Back}(x,\mathbb{X}) = \infty \,\text{ and }\, x \text{ is looping inside } z \oplus\Big[-\frac{1}{2};\frac{1}{2}\Big]^{d} \Big] > 0 ~.
$$
However, on the above event, the random set $Q_{z}(\mathbb{X})$ is infinite which leads to a contradiction with Lemma \ref{lemmetransport}.
\end{dem}

\subsection{An infinite branch of almost looping points}
\label{sect:InfiniteAlmostLooping}

Let us introduce the notion of \textit{almost looping points}. The integer $k$ below is given by \textbf{(LA)}.

\begin{defin}
\label{DefinitionALP}
Let us consider real numbers $0<r<R$, a positive integer $K$, an open ball $A\subset (B(0,r)\times[0;1])^{k}$ and a configuration $\varphi\in\C'$. A marked point $x\in\varphi$ is said a $(r,R,K,A)$-almost looping point of $\varphi$ if:
\begin{itemize}
\item[(i)] $N_{B(x,R)}(\varphi)\le K$;
\item[(ii)] $\forall (x_{1},\ldots,x_{k})\in A_{x}$, we have:
\vskip 0.05cm
\hspace*{1cm} $(ii$-$a)$ $\text{For}(x,\varphi\cup\{x_{1},\ldots,x_{k}\}) \subset \{x,x_{1},\ldots,x_{k}\}$;
\vskip 0.05cm
\hspace*{1cm} $(ii$-$b)$ {There exists a positive constants $C_1$ which does not depend on $x$ such that $\# \text{Back}(x,\varphi\cup\{x_{1},\ldots,x_{k}\}) \geq \# \text{Back}(x,\varphi) - C_{1}$ (in $\mathbb{N}\cup\{\infty\}$).}
\vskip 0.05cm
where $A_{x}=\tau_{\xi}(A)$ with $x=(\xi,\cdot)$.
\end{itemize}
\end{defin}

For a $(r,R,K,A)$-almost looping point, the set $A$ can be interpreted as a suitable region to break the Forward set of $x$ without reducing its Backward set.

The goal of this section is to show Proposition \ref{PropLiggett}. Its proof, given in Section \ref{sect:construction} below, uses intensively \textbf{(LA)} and \textbf{(SA)}.
 
\begin{prop}
\label{PropLiggett}
If $\mathbb{P}(\#\text{For}(\Theta,\mathbb{X}_{\Theta})=\infty)>0$ then, there exists a deterministic quadruplet $(r,R,K,A)$ such that:
\begin{equation}
\label{positive}
\mathbb{P} \Big( \#\{ y \in \text{For}(\Theta,\mathbb{X}_{\Theta}) ;\, y\text{ is a }(r,R,K,A)\text{-almost looping point}\text{ of } \mathbb{X}_{\Theta}\} = \infty \Big) > 0 ~.
\end{equation}
\end{prop}

\subsubsection{Construction of shields}
\label{sect:construction}

Let us first enrich the sequence of events $(\mathscr{E}_{m})_{m\geq 1}$ given by \textbf{(SA)} into a new sequence $(\mathscr{E}'_{m})_{m\geq 1}$. To do it, we need to introduce some definitions.

\begin{defin}
\label{defin:protected}
Let $\varphi\in\C'$ and $m\geq 1$.
\begin{enumerate}
\item A vertex $z\in\Z^{d}$ is said \textit{$m$-shield} for $\varphi$ if $\tau_{-mz}(\varphi)\in\mathscr{E}_{m}$.
\item A vertex $z\in\Z^{d}$ is said \textit{$m$-protecting} for $\varphi$ if for all $y\in\Z^{d}$ such that $\Vert y-z\Vert_{\infty}\in\{0,2\alpha,4\alpha,\ldots,2k\alpha,2(k+1)\alpha,2(k+2)\alpha\}$ (where $\alpha$ is given by \textbf{(SA)}) then $y$ is $m$-shield for $\varphi$.
\item A marked point $x\in\varphi$ is said \textit{$m$-protected} in $\varphi$ if there exists a $m$-protecting vertex $z\in\Z^{d}$ for $\varphi$ such that $x$ belongs to $(mz\oplus[-\alpha m; \alpha m]^{d})\times[0;1]$.
\end{enumerate}
\end{defin}

Roughly speaking, a $m$-protecting vertex $z$ is surrounded by $k+2$ circles (w.r.t. the $\Vert\cdot\Vert_{\infty}$-norm) made up of $m$-shield vertices. Therefore this is also true for a $m$-protected marked point $x$. Thanks to \textbf{(SA)}, each of these circles constitutes an uncrossable obstacle for a single edge $(x,h(\varphi,x))$: see Figure \ref{fig:SA}.

Given a vertex $z\in\mathbb{Z}^{d}$ and two positive integers $m,l\in\mathbb{N}^{*}$, let us set
$$
\text{S}_{m}(z,l) := \left( mz \oplus [-(2l+1)\alpha m;(2l+1)\alpha m]^{d} \right) \times [0;1] ~.
$$
Now, we can define for any integer $m\geq 1$ the event $\mathscr{E}'_{m}$ as the conjunction of the following statements:
\begin{itemize}
\item[$\bullet$] $0$ is $m$-protecting for $\mathbb{X}$ ;
\item[$\bullet$] $\#( \mathbb{X} \cap \text{S}_{m}(0,k+2) ) \le K_{m}$ ;
\item[$\bullet$] $\forall x \in \mathbb{X}_{[-\alpha m;\alpha m]^{d}}$, $\text{rad}(A_{x})>\delta_{m}$ ;
\end{itemize}
where $K_{m},\delta_{m}$ are positive real numbers, and $\text{rad}(A_{x})$ denotes the radius of the open ball $A_{x}$ defined in \textbf{(LA)}.

The use of \textbf{(SA)} allows to assert that the open ball $A_{x}$, for any $m$-protecting point $x\in\varphi$, can be localized from the configuration $\varphi$ only observed through a deterministic and bounded region around $x$. This is the meaning of the next lemma which will be proved in Section \ref{sect:E'mLocal}.

\begin{lem}
\label{lem:E'mLocal}
For any $m\geq 1$, the event $\mathscr{E}'_{m}$ is $\mathscr{S}_{[-\alpha' m;\alpha' m]^{d}}$-measurable where $\alpha':=\alpha (2k+5)$.
\end{lem}

Given a configuration $\varphi\in\mathscr{C}'$ and $x\in\varphi$, we say that $x$ is $m$-\textit{good} for $\varphi$ if there exists a vertex $z\in\mathbb{Z}^{d}$ such that $\tau_{-mz}(\varphi)\in\mathscr{E}'_{m}$ and $x\in\varphi_{mz\oplus[-\alpha m;\alpha m]^{d}}$. Then, for $m$ large enough, the number of $m$-good marked points in a infinite typical forward branch is infinite with positive probability:

\begin{lem}
\label{lem:PercoLiggett}
Assume that $\mathbb{P}(\#\text{For}(\Theta,\mathbb{X}_{\Theta})=\infty)>0$. Then there exists an integer $m_0$ such that, for all $m\ge m_0$,
$$
\mathbb{P} \Big( \#\{x \in \text{For}(\Theta,\mathbb{X}_{\Theta}) ;\, x \text{ is $m$-good for } \mathbb{X}_{\Theta} \} = \infty \Big) > 0 ~.
$$
\end{lem}

\begin{dem}
Let us first prove that the random field
$$
\Xi_{m} := \big( \1_{\{\tau_{-mz}(\mathbb{X})\notin\mathscr{E}'_{m}\}} \big)_{z\in\mathbb{Z}^{d}}
$$
does not percolate in $\mathbb{Z}^{d}$, with probability $1$ and for $m$ large enough. On the one hand, the probability that $0$ is $m$-protecting for $\mathbb{X}$ tends to $1$ as $m\to\infty$ thanks to \textbf{(SA)}. So, one can find a sequence of integers $(K_{m})_{m\ge 1}$ tending to infinity and a sequence of positive real numbers $(\delta_{m})_{m\ge1}$ tending to zero such that the event $\mathscr{E}'_{m}$ also has a probability tending to $1$ as $m\to\infty$. Hence, the probability that $\tau_{-mz}(\mathbb{X})\in\mathscr{E}'_{m} $ goes to $1$ as well. On the other hand, by Lemma \ref{lem:E'mLocal}, $\mathscr{E}'_m$ is $\mathscr{S}_{[-\alpha' m;\alpha' m]^{d}}$-measurable, where $\alpha'=\alpha(2k+5)$. This implies that the event $\{\tau_{-mz}(\mathbb{X})\in\mathscr{E}'_{m}\}$ only depends on the states of vertices $y\in\mathbb{Z}^{d}$ such that $d(y,z)\le 2(2k+5)\alpha$. We can then apply a classic stochastic domination result due to Liggett et al \cite{liggett1997domination}: the random field $\Xi_{m}$ is dominated by an independent Bernoulli field with parameter $p(m)$ tending to $0$ as $m$ tends to infinity. This Bernoulli site percolation, on the lattice $\Z^{d}$ with the $\Vert\cdot\Vert_{\infty}$ graph structure, does not percolate provided $p(m)$ is sufficiently close to $0$. As a consequence, there exists $m_{0}$ such that for all $m\ge m_{0}$,
\begin{equation}
\label{eq:doumayr}
\mathbb{P} \Big(\{z\in\mathbb{Z}^{d}\ ;\ \tau_{-mz}(\mathbb{X})\notin\mathscr{E}'_{m}\}\text{ does not percolate in } \Z^{d}\Big) = 1 ~.
\end{equation}
Combining \eqref{eq:doumayr} with the fact that, thanks to \textbf{(SA)}, it is forbidden to go from one ``bad'' connected component to another one via a single edge, we conclude that the infinite set $\text{For}(\Theta,\mathbb{X}_{\Theta})$ a.s. goes through an infinite number of $m$-protecting vertices $z$ such that $\tau_{-mz}(\mathbb{X}_{\Theta})\in\mathscr{E}'_{m} $.
\end{dem}

This section ends with the proof of Proposition \ref{PropLiggett}.

\begin{dem}[of Proposition \ref{PropLiggett}.]
Assume that $\mathbb{P}(\#\text{For}(\Theta,\mathbb{X}_{\Theta})=\infty)>0$ and let us choose $m\ge m_{0}$ where the integer $m_{0}$ is given by Lemma \ref{lem:PercoLiggett}. Hence, the set of marked points $x\in\text{For}(\Theta,\mathbb{X}_{\Theta})$ which are $m$-goods for $\mathbb{X}_{\Theta}$ is infinite with positive probability. We have to state that an infinite number of them are $(r,R,K,A)$-almost looping points of $\mathbb{X}_{\Theta}$ for some deterministic quadruplet $(r,R,K,A)$.

Let $x\in\text{For}(\Theta,\mathbb{X}_{\Theta})$ be such $m$-good marked point for $\mathbb{X}_{\Theta}$. It is in particular $m$-protected in $\mathbb{X}_{\Theta}$ by some $z\in\Z^{d}$. Lemma \ref{lem:localize} proved in Section \ref{sect:E'mLocal} says that the corresponding set $A_{x}$-- given by \textbf{(LA)} --is included in the marked hypercube $\text{S}_{m}(z,k)^{k}$. So,
$$
\tau_{-\xi}(A_x) \subset \text{S}_{m}(0,k+1)^{k}
$$
where $x=(\xi,\cdot)$.

Now, let us consider a finite covering of the compact set $\text{S}_{m}(0,k+1)^{k} $ by open euclidean balls $\{\mathscr{K}_{j}, 1\le j\le j(m)\}$ of radii $\frac{\delta_{m}}{2}$. Since the open ball $\tau_{-\xi}(A_x)$ is of radius larger than $\delta_m$, it necessarily contains one of the $\mathscr{K}_{j}$'s. We can now conclude by the pigeonhole principle. There exists a deterministic index $1\le j_0\le j(m)$ such that, with positive probability, an infinite number of $m$-good marked points $x\in\text{For}(\Theta,\mathbb{X}_{\Theta})$ satisfy $\tau_{\xi}(\mathscr{K}_{j_{0}}) \subset A_{x}$ where $x=(\xi,\cdot)$. Hence, by \textbf{(LA)}, the deterministic set $\mathscr{K}_{j_{0}}$ satisfies Item $(ii)$ of Definition \ref{DefinitionALP} for all these $m$-good marked points which actually are $(r,R,K,A)$-almost looping points of $\mathbb{X}_{\Theta}$ with $A=\mathscr{K}_{j_{0}}$, $K= K_m$ and any couple $(r,R)$ such that
$$
\mathscr{K}_{j_{0}} \subset (B(0,r)\times[0;1])^{k} \subset (B(0,R)\times[0;1])^{k} \subset \text{S}_{m}(0,k+2)^{k} ~.
$$
\end{dem}



\subsubsection{Proof of Lemma \ref{lem:E'mLocal}}
\label{sect:E'mLocal}

Since $\mathscr{E}_{m}$ is $\mathscr{S}_{[-\alpha m;\alpha m]^{d}}$-measurable by \textbf{(SA)}, the event $\{ 0\text{ is }m\text{-protecting for }\mathbb{X}\}$ is clearly $\mathscr{S}_{[-\alpha' m;\alpha' m]^{d}}$-measurable with $\alpha':=\alpha (2k+5)$. In order to prove that the same holds for $\mathscr{E}'_{m}$, we need to localize the set $A_x$, given by \textbf{(LA)}, for any $x\in\mathbb{X}_{[-\alpha m;\alpha m]^{d}}$.

\begin{lem}
\label{lem:localize}
Let $\varphi\in\C'$ and $m\geq 1$. Let us consider a marked point $x\in\varphi$ which is $m$-protected in $\varphi$ by $0$. Then, the set $A_{x}$ is included in the marked hypercube
$$
A_{x} \subset \text{S}_{m}(0,k)^{k} ~.
$$
\end{lem}

\begin{dem}
Let us consider a $k$-tuple $(x_{1},\ldots,x_{k})$ in $A_{x}$ and assume that one of these vertices does not belong to $\text{S}_{m}(0,k)$. So, one of the $k$ circles of $m$-shield vertices (for $\varphi$) included in $\text{S}_{m}(0,k)$ contains no $x_{i}$'s: there exists $j\in\{1,\ldots,k\}$ such that
\begin{equation}
\label{eq:absx}
\forall i\in\{1,\ldots,k\},\ x_{i} \notin \bigcup_{z'\in V} (mz'\oplus [-\alpha m;\alpha m]^{d}) \times [0;1] ~,
\end{equation}
where $V:=\{z'\in\mathbb{Z}^{d};\ \Vert z'\Vert_{\infty} = 2\alpha j\}$. Thus, let us apply the third item of \textbf{(SA)} to the set $V$  and the configuration $\phi:=\varphi\cup\{x_{1},\dots,x_{k}\}$. We can do it since by hypothesis $\phi$ and $\varphi$ are equal on $mV \oplus [-\alpha m;\alpha m]^{d}$ which means that the vertices of $V$ are still $m$-shield for $\phi$. So, set $\bigcup_{z'\in V}mz'\oplus [-\alpha m;\alpha m]^{d}$ cannot be crossed by an edge of the graph built on $\phi$. But Item $(ii)$ of \textbf{(LA)} guarantees that $x$ and the $x_{i}$'s are in the same connected component in $\phi$ and their Forward sets are included in $\{x,x_{1},\dots,x_{k}\}$. This is possible only if an edge crosses the $\bigcup_{z'\in V}mz'\oplus [-\alpha m;\alpha m]^{d}$. We just have seen that such situation was forbidden.
\end{dem}

{\begin{lem}
\label{lem:localize2}
Let $\varphi\in\mathscr{C}'$ such that the vertex $0$ is $m$-protecting. Let $x\in\varphi_{[-\alpha m;\alpha m]^{d}}$. Then, for any $k$-tuple $(x_{1},\dots,x_{k})\in \text{S}_{m}(0,k)^{k}$, the control of the three items of \textbf{(LA)} only depends on $\phi\cap\text{S}_{m}(0,k+2)$ where $\phi:=\varphi\cup\{x_{1},\dots,x_{k}\}$.
\end{lem}}

\begin{figure}[!ht]
\begin{center}
\vspace{0.5cm}
\psfrag{x}{\small{$x$}}
\psfrag{x1}{\small{$x_{1}$}}
\psfrag{x2}{\small{$x_{2}$}}
\psfrag{x3}{\small{$x_{3}$}}
\psfrag{Shield1}{\small{$\text{S}_{m}(0,k)$}}
\psfrag{Shield2}{\small{$\text{S}_{m}(0,k+2)^{c}$}}
\includegraphics[width=10cm, height=10cm]{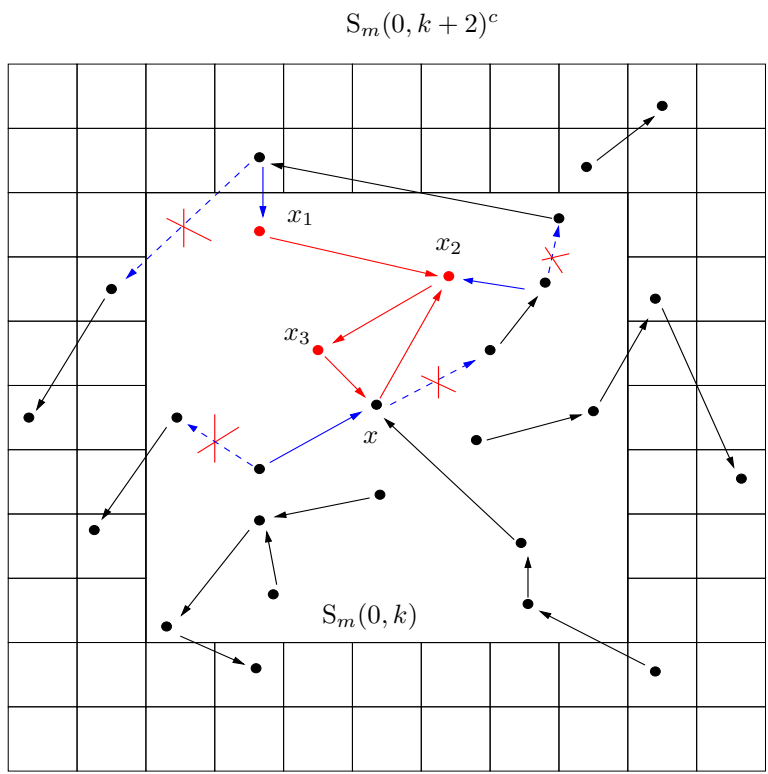}
\caption{\label{fig:shieldloop} On this picture, $k=3$ marked points, denoted by $x_{1}$, $x_{2}$ and $x_{3}$, are added to the configuration $\varphi$ in which $x$ is $m$-protected by $0$. Remark that the three items of \textbf{(LA)} are satisfied with the new configuration $\phi:=\varphi\cup\{x_{1},\dots,x_{k}\}$. The Forward sets of $x$ and the $x_i$'s are included in $\{x,x_{1},x_{2},x_{3}\}$. Moreover, the add of the $x_i$'s destroys some edges started from marked points of $\varphi$ and also creates new ones. This second phenomenon explains that the Backward set of $x$ increases here-- which is allowed by the third item of \textbf{(LA)}.}
\end{center}
\end{figure}

Combining the two previous lemmas, we deduce Lemma \ref{lem:E'mLocal}. This section ends with the proof of Lemma \ref{lem:localize2}.

\begin{dem}
Lemma \ref{lem:localize} ensures that the set $A_{x}$ is included in the marked hypercube $\text{S}_{m}(0,k)^{k}$. This implies the following important fact. The vertices of $V:=\{z'\in\mathbb{Z}^{d}\ ;\ \|z'\|_{\infty}=2\alpha (k+1)\}$ which were $m$-shield for $\varphi$ (since $0$ is $m$-protecting for $\varphi$) are still $m$-shield for $\phi$. Indeed, the configurations $\varphi$ and $\phi$ coincide on the set $mV \oplus [-\alpha m;\alpha m]^{d}$ and $\mathscr{E}_{m}$ is $\mathscr{S}_{[-\alpha m;\alpha m]^{d}}$-measurable.

Let us pick $(x_{1},\dots,x_{k})\in \text{S}_{m}(0,k)^{k}$. Look at Figure \ref{fig:shieldloop} for an example. Let us first check that $\text{For}(x,\phi)\subset\{x,x_{1},\dots,x_{k}\}$. It is sufficient to determine the oriented edges starting from the points $x,x_{1},\dots,x_{k}$ in the graph built on $\phi$. By the previous remark, we can use the third item of \textbf{(SA)} w.r.t. the set $V$: the set $mV \oplus [-\alpha m;\alpha m]^{d}$ splits the space into two disjoint connected components where the bounded one is $\text{S}_{m}(0,k)$. The oriented edges starting from $\phi\cap\text{S}_{m}(0,k)$ do not depend on $\phi\cap\text{S}_{m}(0,k+1)^{c}$. Hence, we can determine the outgoing edges of $x,x_{1},\dots,x_{k}$ without consider $\phi$ outside $\text{S}_{m}(0,k+1)$. So, the control of Item $(i)$ of \textbf{(LA)} only depends on $\phi\cap \text{S}_{m}(0,k+1)$. {The same argument works as well for Item $(ii)$ of \textbf{(LA)} and leads to the same conclusion.}

{It then remains to prove that the verification of Item $(iii)$ of \textbf{(LA)} only depends on $\phi\cap\text{S}_{m}(0,k+2)$. To do it, we have to prove that the outgoing vertex $h(\phi,y)$, for any vertex $y$ such that $h(\phi,y)\not=h(\varphi,y)$ i.e. for any $y$ whose outgoing edge is altered when adding $x_{1},\dots,x_{k}$, can be identified thanks to $\phi\cap\text{S}_{m}(0,k+2)$. Let us pick such a vertex $y$. The previous argument involving vertices of the set $V=\{z'\in\mathbb{Z}^{d}\ ;\ \|z'\|_{\infty}=2\alpha (k+1)\}$  (combined with \textbf{(SA)}) forces $y$ to belong to $\text{S}_{m}(0,k+1)$. Thus, the same argument used this time with the set $V':=\{z'\in\mathbb{Z}^{d}\ ;\ \|z'\|_{\infty}=2\alpha (k+2)\}$ (this is the reason why we need $k+2$ circles of $m$-shield vertices in the Definition \ref{defin:protected}) implies that the outgoing vertex $h(\phi,y)$ does not depend on what happens outside $\text{S}_{m}(0,k+2)$.}
\end{dem}

\subsection{From Forward set to Backward set} 
\label{sect:ForwardToBackward}

In this section, it is stated that the mean size of the Backward set of a typical almost looping point is infinite whenever the Forward set of a typical marked point contains an infinite number of almost looping points with positive probability. Above all, this result allows to convert a forward result to a backward one.

\begin{prop}
\label{prop:contradiction}
If there exist a quadruplet $(r,R,K,A)$ such that
\begin{equation}
\label{equa}
\mathbb{P} \Big( \#\{ y \in \text{For}(\Theta,\mathbb{X}_{\Theta}) ;\, y\text{ is a $(r,R,K,A)$-almost looping point of } \mathbb{X}_{\Theta} \} = \infty \Big) > 0 ~,
\end{equation}
then
\begin{equation}
\label{Equa}
\mathbb{E} \left[ \#\text{Back}(\Theta,\mathbb{X}_{\Theta}) \1_{\{\Theta\text{ is a $(r,R,K,A)$-almost looping point of } \mathbb{X}_{\Theta} \}} \right] = \infty ~. 
\end{equation}
\end{prop}

\begin{dem} 
Let us fix parameters $r,R,K,A$ such that (\ref{equa}) holds. Let us denote by $\text{For}^{\ast}(\Theta,\mathbb{X}_{\Theta})$ the subset of $\text{For}(\Theta,\mathbb{X}_{\Theta})$ made up of $(r,R,K,A)$-almost looping points:
$$
\text{For}^{\ast}(\Theta,\mathbb{X}_{\Theta}) = \{ y \in \text{For}(\Theta,\mathbb{X}_{\Theta}) ;\, y\text{ is a $(r,R,K,A)$-almost looping point of }\mathbb{X}_{\Theta} \} ~.
$$
We need to bound from below the density of $\text{For}^{\ast}(\Theta,\mathbb{X}_{\Theta})$. The next result will be proved at the end of the section.

\begin{lem} 
\label{lemmadense}
There exists a function $g:\mathbb{N}\to\mathbb{N}$ such that $\lim_{n\rightarrow \infty} g(n)=\infty$ and
$$
\mathbb{P} \Big( \forall n\ge 1, \, \#\Big( \text{For}^{\ast}(\Theta,\mathbb{X}_{\Theta}) \cap ([-n;n]^{d}\times[0;1])\Big) \ge g(n) \Big) > 0 ~.
$$
\end{lem}

We will say that a marked point $x$ is \textit{dense} for $\mathbb{X}$ if for any $n\geq 1$,
$$
\#\Big( \text{For}^{\ast}(x,\mathbb{X}) \cap ([-n;n]^{d}\times[0;1])\Big) \ge g(n) ~.
$$
Lemma \ref{lemmadense} asserts that, with positive probability, $\Theta$ is dense in $\mathbb{X}_{\Theta}$.

Let $\Lambda_{n}=[-n;n]^{d}\times[0;1]$ for $n\geq 1$. The {Campbell Mecke Formula (see for example \cite{LastPenroseBook})} allows to write:
$$
\mathbb{E} \left[\#\text{Back}(\Theta,\mathbb{X}_{\Theta}) \1_{\{ \Theta\text{ is a $(r,R,K,A)$-almost looping point of }\mathbb{X}_{\Theta} \}} \right] \hspace*{4cm}
$$
$$
\hspace*{2cm}= \frac{1}{(4n)^{d}}\mathbb{E} \left[ \sum_{x\in \mathbb{X}\cap \Lambda_{2n}} \#\text{Back}(x,\mathbb{X}) \1_{\{ x\text{ is a $(r,R,K,A)$-almost looping point of } \mathbb{X} \}} \right] ~.
$$
If the marked point $x\in\Lambda_{n}$ is dense for $\mathbb{X}$ then there exist at least $g(n)$ marked points in $\Lambda_{2n}\cap\mathbb{X}$ having $x$ in their Backward sets. This is here that the passage from the Forward set to the Backward set happens. Henceforth,
\begin{eqnarray*}
\mathbb{E} \left[\#\text{Back}(\Theta,\mathbb{X}_{\Theta}) \1_{\{ \Theta\text{ is a $(r,R,K,A)$-almost looping point of }\mathbb{X}_{\Theta} \}} \right] & \geq & \frac{g(n)}{(4n)^{d}} \mathbb{E} \left[ \sum_{x\in \Lambda_{n}\cap\mathbb{X}} \1_{\{ x \text{ is dense for }\mathbb{X} \}} \right]\\
& = & \frac{1}{2^{d}} g(n) \mathbb{P}(\Theta \text{ is dense in } \mathbb{X}_{\Theta}).
\end{eqnarray*}
Let us tend $n$ to infinity. By Lemma \ref{lemmadense}, (\ref{Equa}) follows.
\end{dem}

\begin{dem}[of Lemma \ref{lemmadense}]
Let $(a_{n})_{n \ge 1}$ be a sequence of positive real numbers whose sum $\sum_n a_n$ equals $\alpha/2$ where $\alpha=\mathbb{P}(\#\text{For}^{\ast}(\Theta,\mathbb{X}_{\Theta})=\infty)>0$ by hypothesis. Thus, we define a (nondecreasing) sequence of integers $(n_{k})_{k\ge 1}$ and a sequence of events $(B_{k})_{k\ge 0}$ as follows. At first, let us define $B_{0}$ as the event $\{\#\text{For}^{\ast}(\Theta,\mathbb{X}_{\Theta})=\infty\}$ and $n_{1}$ as the first integer $n$ such that
$$
\mathbb{P} \Big( B_{0} \cap \Big\{ \# \Big( \text{For}^{\ast}(\Theta,\mathbb{X}_{\Theta}) \cap ([-n;n]^{d}\times[0;1]) \Big) \ge 1 \Big\} \Big) \ge \alpha-a_{1} ~.
$$
Since the above probability tends to $\alpha$ as $n$ tends to infinity, $n_{1}$ is well defined. We also denote by $B_{1}$ the following event:
$$
B_{1} = B_{0} \cap \Big\{ \#\Big( \text{For}^{\ast}(\Theta,\mathbb{X}_{\Theta}) \cap ([-n_{1};n_{1}]^{d}\times[0;1]) \Big) \ge 1 \Big\} ~.
$$
Thus, for any $k\geq 2$, we define by induction the integer $n_{k}$ as the first integer $n$ such that
$$
\mathbb{P} \Big( B_{k-1} \cap \Big\{ \# \Big( \text{For}^{\ast}(\Theta,\mathbb{X}_{\Theta}) \cap ([-n;n]^{d}\times[0;1]) \Big) \ge k\Big\} \Big) \ge \alpha-\sum_{1\leq i\leq k-1} a_{i} ~.
$$
We also set
$$
B_{k} = B_{k-1} \cap \Big\{ \#\Big( \text{For}^{\ast}(\Theta,\mathbb{X}_{\Theta}) \cap ([-n_{k};n_{k}]^{d}\times[0;1]) \Big) \ge k \Big\} ~.
$$
Finally, let us define a function $g$ by $g(n)=\sum_k \1_{[n_{k}, n_{k+1}[}$ $ (n)$. The value $g(n)$ gives the number of $n_{k}$'s smaller than $n$. The function $g$ has been built so as to satisfy
$$
\mathbb{P} \left( \bigcap_{k \ge 1} B_{k} \right) \leq \mathbb{P} \left( \forall n\ge 1,\ \#\text{For}^{\ast}(\Theta,\mathbb{X}_{\Theta})\cap ([-n;n]^{d}\times[0;1])\ge g(n) \right) ~.
$$
Finally, to conclude, it suffices to remark that $\cap_{k \ge 1} B_{k}$ occurs with probability larger than $\alpha-\sum_{k\ge 1} a_{k}$ which is equal to $\frac{\alpha}{2}$.
\end{dem}

\subsection{From almost looping points to looping points}
\label{sect:FromAlmostLoopingToLooping}

A \textit{looping point} is a marked point whose Forward set admits a localized loop-- actually, only its center of mass is localized.

\begin{defin}
\label{defi:LoopingPoint}
Let $r<R$ be some positive real numbers and $K$ be a positive integer. Let $\varphi\in\C'$. A marked point $x\in\varphi$ is a $(r,R,K)$-looping point of $\varphi$ if
\begin{itemize}
\item[(i)] $N_{B(x,R)}(\varphi)\le K$;
\item[(ii)] $x$ is looping inside the ball $B(x,r)$ for $\varphi$ (definition given in Section \ref{sect:ForwardtoBackward}).
\end{itemize}
\end{defin}

So, a $(r,R,K,A)$-almost looping point $x$ of $\varphi$ becomes a $(r,R,K+k)$-looping point of $\varphi\cup\{x_{1},\ldots,x_{k}\}$ when the $k$ marked points $x_{1},\ldots,x_{k}$ are added in $A_x$. Proposition \ref{prop:FromAlmostLoopingToLooping} establishes a link between almost looping points and looping points.

\begin{prop}
\label{prop:FromAlmostLoopingToLooping}
If there exist a quadruplet $(r,R,K,A)$ such that
\begin{equation}
\label{EALP}
\mathbb{E} \left[ \#\text{Back}(\Theta,\mathbb{X}_{\Theta})\1_{\lbrace \Theta\text{ is a }(r,R,K,A)-\text{almost looping point for } \mathbb{X}_{\Theta} \rbrace} \right] = \infty
\end{equation}
then
\begin{equation}
\label{ELP}
\mathbb{E} \left[ \#\text{Back}(\Theta,\mathbb{X}_{\Theta})\1_{\lbrace \Theta\text{ is a } (r,R,K+k)-\text{looping point for } \mathbb{X}_{\Theta}\rbrace} \right] = \infty ~.
\end{equation}
\end{prop}

The rest of this section is devoted to the proof of Proposition \ref{prop:FromAlmostLoopingToLooping}. Let us first state without proof a technical lemma (Exercise 4.10 of \cite{LastPenroseBook}).

\begin{lem}
\label{LemDensity}
Let $\Lambda$ be a bounded Borelian of $\Rd$ with positive Lebesgue measure. Let us denote by $U$ the uniform distribution on $\Lambda$. Let us consider $(X_{i})_{1\le i\le k}$ i.i.d. random vectors on $\Lambda\times[0;1]$ with distribution $U\otimes Q$ which are also independent with $\mathbb{X}$. Let us set $\mathbb{X}'=\mathbb{X}\cup\{X_{1},\ldots,X_{k}\}$ the extended point process. Then the law $\Pi'$ of $\mathbb{X}'$ is absolutely continuous with respect to the law $\Pi$ of $\mathbb{X}$ (i.e. the Poisson Point distribution) with density
\begin{equation}
\label{density}
\frac{\Pi^{'}(d\varphi)}{\Pi(d\varphi)} = \frac{1}{\lambda_{d}(\Lambda)^{k}} N_{\Lambda}(\varphi) \left( N_{\Lambda}(\varphi)-1 \right) \ldots \left( N_{\Lambda}(\varphi)-k+1 \right) ~.
\end{equation}
\end{lem}

\begin{dem}[of Proposition \ref{prop:FromAlmostLoopingToLooping}.]
Let a quadruplet $(r,R,K,A)$ such that (\ref{EALP}) occurs. Our goal is to prove that the expectation in (\ref{ELP}), denoted by $\mathscr{I}$, is infinite.
Lemma \ref{LemDensity} applied to the set $\Lambda=B(0,R)$ allows to write:
\begin{eqnarray*}
\mathscr{I} & \ge & \mathbb{E}\left[\#\text{Back}(\Theta,\mathbb{X}_{\Theta})\1_{\lbrace \text{$\Theta$ is looping inside $B(0,r)$ for $\mathbb{X}_{\Theta}$} \rbrace} \1_{\lbrace k\le N_{B(0,R)}(\mathbb{X})\le K+k-1\rbrace} \right] \\
& = & \int_{\C'} \#\text{Back}(\Theta,\varphi_{\Theta}) \1_{\lbrace \text{$\Theta$ is looping inside $B(0,r)$ for $\varphi_{\Theta}$} \rbrace} \1_{\lbrace k \le N_{B(0,R)}(\varphi)\le K+k-1\rbrace} \Pi(d\varphi) \\
& = & \int_{\C'} \#\text{Back}(\Theta,\varphi_{\Theta}) \1_{\lbrace \text{$\Theta$ is looping inside $B(0,r)$ for $\varphi_{\Theta}$} \rbrace} \1_{\lbrace k\le N_{B(0,R)}(\varphi)\le K+k-1\rbrace} \frac{\Pi^{'}(d\varphi)}{f(\varphi)} ~,
\end{eqnarray*}
where $f(\varphi)$ is the density given in (\ref{density}). Provided $k\le N_{B(0,R)}(\mathbb{X})\le K+k-1$, the ratio $1/f(\varphi)$ is larger than some constant $C>0$. It follows:\begin{eqnarray*}
\mathscr{I} & \ge & C \int_{\C'} \#\text{Back}(\Theta,\varphi_{\Theta}) \1_{\lbrace \text{$\Theta$ is looping inside $B(0,r)$ for $\varphi_{\Theta}$} \rbrace} \1_{\lbrace k\le N_{B(0,R)}(\varphi)\le K+k-1\rbrace} \Pi'(d\varphi) \\
& \ge & C \int_{A} \mathbb{E} \Big[ \#\text{Back}(\Theta,\mathbb{X}_{\Theta}\cup\{x_{1},\ldots,x_{k}\}) \\
& & \qquad \times \1_{\lbrace \text{$\Theta$ is looping inside $B(0,r)$ for $\mathbb{X}_{\Theta}\cup\{x_{1},\ldots,x_{k}\}$} \rbrace} \1_{\lbrace N_{B(0,R)}(\mathbb{X})\le K-1\rbrace} \Big] (U\otimes Q)^{k} [dx_{1}\ldots dx_{k}]\\
& \ge & {C (U\otimes Q)^{k}(A) \left( \mathbb{E} \left[ \#\text{Back}(\Theta,\mathbb{X}_{\Theta}) \1_{\lbrace \Theta \text{ is a } (r,R,K,A)\text{-almost looping point}\text{ of } \mathbb{X}_{\Theta} \rbrace} \right] \, - \, C_1 \right)}
\end{eqnarray*}
which is infinite by hypothesis (we have $(U\otimes Q)^{k}(A)>0$). This concludes the proof. It is worth pointing out here that the condition $(ii$-$b)$ is used to obtain the latter inequality.
\end{dem}

\subsection{Proof of Theorem \ref{mainTH}: conclusion.}
\label{sect:LoopingPoint}

Another use of the mass transport principle, especially Lemma \ref{lemmetransport}, leads to the next result: the Backward set of a typical looping point has a finite mean size.

\begin{prop}
\label{PropBLP}
Any triplet $(r,R,K)$ satisfies
\begin{equation}
\label{BLP}
\mathbb{E} \Big[ \#\text{Back}(\Theta,\mathbb{X}_{\Theta}) \1_{ \{\Theta \text{ is a $(r,R,K)$-looping point of } \mathbb{X}_{\Theta}\}} \Big] < \infty ~.
\end{equation}
\end{prop}

\begin{dem}
Let $r<R$ be some positive real numbers and $K$ be a positive integer. Let us pick $\epsilon>0$ small enough so that $D_{1}\subset D_{2}$ where
$$
D_{1} = \bigcup_{\eta\in\big[-\frac{\epsilon}{2};\frac{\epsilon}{2}\big]^{d}} B(\eta,r) \; \text{ and } \; D_{2} = \bigcap_{\eta\in\big[-\frac{\epsilon}{2};\frac{\epsilon}{2}\big]^{d}} B(\eta,R) ~.
$$
Let $\mathscr{I}$ be the expectation in (\ref{BLP}). Using the Campbell Mecke Formula on the set $M=[-\frac{\epsilon}{2};\frac{\epsilon}{2}]^{d}\times[0;1]$, we can write:
\begin{eqnarray}
\label{BackLoopInside}
\mathscr{I} & = & \frac{1}{\epsilon^{d}} \mathbb{E} \left[ \sum_{x\in \mathbb{X}\cap M} \#\text{Back}(x,\mathbb{X}) \1_{\{\text{For}(x,\mathbb{X}) \text{ is looping inside } B(x,r) \}} \1_{\{ N_{B(x,R)}(\mathbb{X})\le K \}}\right] \nonumber \\
& \le & \frac{1}{\epsilon^{d}} \mathbb{E} \left[ \sum_{x\in \mathbb{X}\cap M}  \#\text{Back}(x,\mathbb{X}) \1_{\{ \text{For}(x,\mathbb{X})\text{ is looping inside } D_{1} \}} \1_{\{ N_{D_{2}}(\mathbb{X})\le K \}} \right] \\
& \le & \frac{K}{\epsilon^{d}} \mathbb{E} \Big[  \#\{ y \in \mathbb{X} ;\, \text{For}(y,\mathbb{X}) \text{ is looping inside } D_ {1} \} \Big] \nonumber~,
\end{eqnarray}
since each marked point $y$ whose Forward set is looping inside $D_ {1}$ is counting at most $K$ times in the sum of (\ref{BackLoopInside}). Thus, Lemma \ref{lemmetransport} allows to conclude.
\end{dem}

It only remains to combine the different pieces of the proof of Theorem \ref{mainTH}.

{\begin{dem}[of Theorem \ref{mainTH}.]
By Proposition \ref{PropTM}, it is enough to show that a.s. the Forward set $\text{For}(x,\mathbb{X})$ of any $x\in\mathbb{X}$ is finite. The {Campbell Mecke Formula} gives:
$$
\mathbb{E} \left[ \sum_{x\in\mathbb{X}}  \1_{\{ \#\text{For}(x,\mathbb{X})=\infty \}} \right] = \int \mathbb{P} \big( \#\text{For}(x,\mathbb{X}\cup\{x\})=\infty \big) \, \lambda_{d} \otimes Q (dx) ~.
$$
By stationarity, we have to prove that $\mathbb{P}(\#\text{For}(\Theta,\mathbb{X}_{\Theta})=\infty)=0$ where $\Theta=(0,U)$ and $U$ is a uniform r.v. in $[0;1]$, and $\mathbb{X}_{\Theta}=\mathbb{X}\cup\{\Theta\}$.

If the probability $\mathbb{P}(\#\text{For}(\Theta,\mathbb{X}_{\Theta})=\infty)$ was positive then, by Propositions \ref{PropLiggett}, \ref{prop:contradiction} and \ref{prop:FromAlmostLoopingToLooping}, there would exist a triplet $(r,R,K)$ such that
$$
\mathbb{E} \left[ \#\text{Back}(\Theta,\mathbb{X}_{\Theta})\1_{\lbrace \Theta\text{ is a } (r,R,K+k)-\text{looping point for } \mathbb{X}_{\Theta}\rbrace} \right] = \infty
$$
which is in contradiction with Proposition \ref{PropBLP}.
\end{dem}

\section{Proof of Theorem \ref{theo:models}}
\label{SectionPC}

This section is devoted to the verifications of \textbf{(LA)} and \textbf{(SA)} for the \textbf{Line-segment model} and the \textbf{Navigation model}.

\subsection{Line-segment model}
\label{sect:HypoGrowing}

Let us introduce some notations decribing the geometry of this model. Given a marked point $x=(\xi,u)\in\mathbb{R}^{2}\times[0;1]$, we denote by $l(x)=\{ \xi+t\overrightarrow{u}, t\in\mathbb{R}_{+}\}$ the (semi-infinite) ray starting from $\xi$ in the direction $\overrightarrow{u}=(\cos(2\pi u),\sin(2\pi u))$. Thus, let us set
$$
l(\mathbb{X}) = \bigcup _{x\in\mathbb{X}} l(x) ~.
$$
For $\varphi\in\mathscr{C}^{'}$ and $x=(\xi,u)\in\varphi$, we denote by $h_{g}(\varphi,x)\in\mathbb{R}^{2}$ the intersection point between $l(x)$ and $l(h(\varphi,x))$. Roughly speaking, $h_{g}(\varphi,x)$ represents the impact point of the stopped segment starting from $\xi$ on the stopping segment starting from $\xi'$, where $(\xi',\cdot)=h(\varphi,x)$.

\subsubsection{Loop assumption}

Let us prove that the Line-segment model satisfies \textbf{(LA)} with $k=3$.

\begin{prop}
\label{prop:CyclingAssGrowing}
With probability $1$, $\mathbb{X}$ is a $3$-looping configuration.
\end{prop}

Consider a configuration $\varphi\in\mathscr{C}'$ and an element $x=(\xi,u)\in\varphi$. The first step consists in stating that only a finite number of growing segments are stopped by $[\xi;h_{g}(x,\varphi)]$. See Lemma \ref{lem:DefLine} below. This will allow us to exhibit a small region close to the impact point $h_{g}(\varphi,x)$ where we could easily add a loop made up of $3$ segments, that the growing segment $x$ will hit. See Figure \ref{fig:cicleass}.

\begin{figure}[!ht]
\begin{center}
\psfrag{x}{\small{$x=(\xi,u)$}}
\psfrag{y}{\small{$h(\varphi,x)$}}
\psfrag{z1}{\small{$z_{1}$}}
\psfrag{z2}{\small{$z_{2}$}}
\includegraphics[width=9cm,height=6cm]{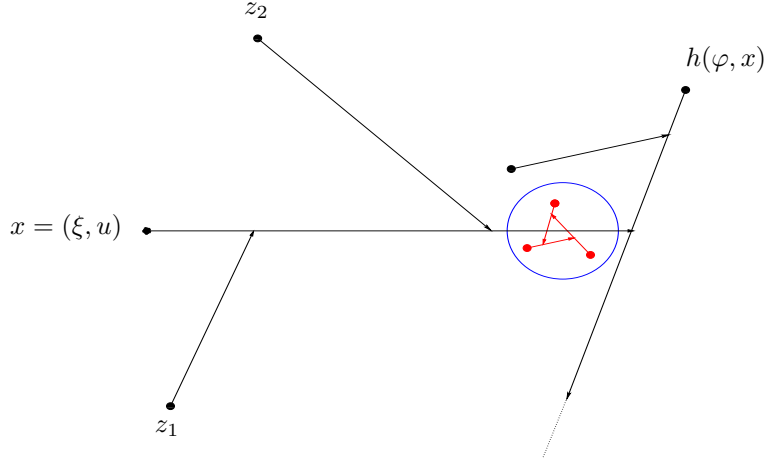}
\end{center}
\caption{\label{fig:cicleass} On this picture, the set $\text{Back}^{-1}(\varphi,x)=\{y\in\varphi, h(\varphi,y)=x\}$ is made up of the marked points $z_{1}$ and $z_{2}$. The blue circle delimits the ball $B(w,r')$. It is worth pointing out here that $B(w,r')$ has to avoid any growing segment $y\in\mathbb{X}$, and not only $y\in\text{Back}^{-1}(\varphi,x)$, so that the added marked points (in red) may stop the growing segment $x$ without decreasing its Backward set.}
\end{figure}

Let us denote by $\mathscr{B}$ the set of open discs in $\mathbb{R}^{2}$, {and by $\mathscr{C}''$ the following measurable set:
$$
\mathscr{C}'' =  \Big\{ \varphi\in\mathscr{C}'\ ;\ \forall B\in\mathscr{B},\ \#\{x=(\xi,u)\in\varphi\ :\ [\xi;h_{g}(\varphi,x)]\cap B \not=\emptyset\} < +\infty \Big\} ~.
$$
Remark that the measurability of $\mathscr{C}''$ is based on the one of $h_{g}:(\varphi,x)\in\mathscr{C}'\times\mathbb{R}^{2}\to h_{g}(\varphi,x)\in\mathbb{R}^{2}$. For this, we refer the reader to Section 4 of \cite{daley2014two}.}

\begin{lem}
\label{lem:DefLine}
With probability $1$, $\mathbb{X}$ belongs to $\mathscr{C}''$.
\end{lem}

We are now able to prove Proposition \ref{prop:CyclingAssGrowing}.

\begin{dem}
Let us first introduce the set of marked points of a configuration $\varphi$ which are stopped by $x$:
$$
\text{Back}^{-1}(\varphi,x) = \{ y\in\varphi , \, h(\varphi,y)=x \} ~.
$$
Let $\varphi\in\mathscr{C}''$ and $x=(\xi,u)\in\varphi$. Recall that the impact point $h_{g}(\varphi,y)$ of a marked point $y\in\text{Back}^{-1}(\varphi,x)$ belongs to $[\xi;h_{g}(\varphi,x)]$. Since $\text{Back}^{-1}(\varphi,x)$ is finite, we can exhibit a positive real number $r$ such that
\begin{equation}
\label{EmptyZone}
[h_{g}(\varphi,x)-r\overrightarrow{u};h_{g}(\varphi,x)] \cap \left\{ h_{g}(\varphi,y) , \, y\in\text{Back}^{-1}(\varphi,x) \right\} = \emptyset ~.
\end{equation}
Let us set $w=h_{g}(\varphi,x)-\frac{r}{2}\overrightarrow{u}$. Statement (\ref{EmptyZone}) ensures the existence of a positive radius $r'$ such that
\begin{equation}
\label{AvoidCond}
\left(\bigcup_{y=(\eta,v)\in\varphi}  [\eta;h_{g}(\varphi,y)]\right)\cap B(w,r') = [\xi;h_{g}(\varphi,x)]\cap B(w,r')
\end{equation}
on the one hand, and
\begin{equation}
\label{NoStopByX}
2r' \le \|\xi - w\|_{2} - r'
\end{equation}
on the other hand. The obtained ball $B(w,r')$ is actually a suitable region in which we could create an obstacle for the growing segment $x$ without altering any other growing segment. Besides, condition (\ref{NoStopByX}) ensures that any marked point added in $B(w,r')$ could not be stopped by the growing segment $x$.

Let us consider the set $\mathscr{A}_{x}$ of triplets  $(x_{0},x_{1},x_{2})\in(B(w,r')\times [0;1])^{3}$ such that:
\begin{itemize}
\item[$(i)$] $h(x_{i},\varphi\cup\{x_{0},x_{1},x_{2}\}) = x_{i+1}$ for $i=0,1,2$ (where the index $i+1$ is taken modulo $3$);
\item[$(ii)$] The triangle defined by the vertices $h_{g}(x_{i},\varphi\cup\{x_{0},x_{1},x_{2}\})$, $i=0,1,2$, is included in $B(w,r')$ and contains the center $w$.
\end{itemize}

It is not difficult too see that $\mathscr{A}_{x}$ contains a non-empty open set $A_{x}\subset (B(w,r^{'})\times[0;1])^{3}$. When a triplet $(x_{0},x_{1},x_{2})\in A_{x}$ is added to $\varphi$, then by (\ref{NoStopByX}), $(i)$ and $(ii)$, the growing segment $x$ hits the loop produced by $x_{0},x_{1},x_{2}$:
\begin{eqnarray*}
\text{For}(x,\varphi\cup\{x_{0},x_{1},x_{2}\}) &=& \{x,x_{0},x_{1},x_{2}\},\\
\forall 0\le i\le 2,\ \text{For}(x_{i},\varphi\cup\{x_{0},x_{1},x_{2}\}) &=& \{x_{0},x_{1},x_{2}\}.
\end{eqnarray*}

Then the first two items of \textbf{(LA)} are checked. Morever, condition (\ref{AvoidCond}) in conjunction with $(i)$ and $(ii)$ imply that no growing segment except $x$ is changing by the adding marked points $\{x_{0},x_{1},x_{2}\}$, the third item of \textbf{(LA)} is also checked and:
$$
\text{Back}(x,\varphi\cup\{x_{0},x_{1},x_{2}\}) = \text{Back}(x,\varphi) ~.
$$
This achieves the proof of Proposition \ref{prop:CyclingAssGrowing}.
\end{dem}

\begin{dem}[of Lemma \ref{lem:DefLine}.]
Using classical arguments, it is sufficient to prove that:
$$
\mathbb{P}\left( \#\{x=(\xi,\varpi)\in\mathbb{X}\ ;\ [\xi;h_{g}(\mathbb{X},x)]\cap B(0,1)\not=\emptyset\}<+\infty\right)=1 ~.
$$
We will show that:
$$
\mathscr{E}=\mathbb{E}\left(\#\{x=(\xi,u)\in\mathbb{X}\ ;\ [\xi;h_{g}(\mathbb{X},x)]\cap B(0,1)\not=\emptyset\}\right)<+\infty ~.
$$
Let us apply the Campbell-Mecke formula:
\begin{eqnarray*}
\mathscr{E} &=& z\pi+z\int_{(B(0,1)\times[0;1])^{c}}\mathbb{E}
\left(\1_{\{\left[x;,h_{g}(\mathbb{X}\cup\{x\},x)\right]\cap B(0,1)\not=\emptyset\}}\right)\lambda_{2}(d\xi)Q(du),\\
& \le & z\pi+2\pi z \int^{1}_{0}\left(\int_{1}^{+\infty} \mathbb{P}\left(\|\xi-h_{g}(\mathbb{X}\cup\{x\},x)\|\ge r-1\right)rdr\right)Q(du),\\
& \le & z\pi+2\pi z\int^{1}_{0}\left(\int_{1}^{+\infty} \mathbb{P}\left(\|h_{g}\left(\mathbb{X}\cup\{(0,u)\},(0,u)\right)\|\ge r-1\right)rdr\right)Q(du).
\end{eqnarray*} 

By isotropy, for all $u\in[0,1]$ and for all $r>0$,
$$
\mathbb{P} \left(\|h_{g}\left(\mathbb{X}\cup\{(0,u)\},(0,u)\right)\|\ge r\right) = \mathbb{P} \left(\|h_{g}(\mathbb{X}\cup\{(0,0)\},(0,0))\|\ge r \right) ~,
$$
where $(0,0)$ denotes the marked vertex located at the origin with direction $(1,0)$. Schreiber \& Soja have proved (Theorem 4 in \cite{schreiber2010limit}) that there exist $c,c'>0$ such that the probability $\mathbb{P}\left(\|h_{g}(\mathbb{X}\cup\{(0,0)\},(0,0))\|\ge r\right)$ is smaller than $c\text{e}^{-c'r} $ for all $r\ge 0$. This exponential decay ensures that $\mathscr{E}$ is finite.
\end{dem}

\subsubsection{Shield assumption}

In order to prove \textbf{(SA)}, we need to construct building blocks called \textit{shield hexagons} which together will formed uncrossable walls. To do it, let us start with introducing an hexagonal tessellation.

Let us consider the triangular lattice whose vertex set is
$$
\Pi = \left\{ a\overrightarrow{i}+b\overrightarrow{j}\ :\ a,b\in\mathbb{Z} \right\} ~,
$$
where $\overrightarrow{i}=(\sqrt{3}.\cos(\frac{\pi}{6}),\sqrt{3}.\sin(\frac{\pi}{6}))$ and $\overrightarrow{j}=(0,\sqrt{3})$. The usual graph distance on $\Pi$ is denoted by $d_{\Pi}$. We also denote resp. by $B^{n}(z)$ and $S^{n}(z)$ the (closed) ball and sphere with center $z$ and radius $n$ w.r.t. $d_{\Pi}$.

For any $z\in\Pi$, let $\H(z)$ be the Voronoi cell of $z$ w.r.t. the vertex set $\Pi$:
$$
\H(z) = \Big\{ y\in\mathbb{R}^{2},\ \|y-z\|_{2} \le \inf_{w\in\Pi\setminus\{z\}} \|y-w\|_{2}\Big\} ~.
$$
The set $\H(z)$ is a regular hexagon centred at $z$. For any integer $n>0$, let us introduce the hexagonal complex of size $n$ centred in $z$ as
$$
\H^{n}(z) = \bigcup_{y\in B^{n}(z)} \H(y) ~.
$$
Given $\xi\in\mathbb{R}^{2}$, we also set $\H^{n}(\xi)=\H^{n}(0)+\xi$. Finally, for any integer $n>0$, we define the hexagonal ring $C_{n}(\xi)$ by $C_{n}(\xi)=\H^{n}(\xi)\setminus \H^{n-1}(\xi)$ (with $\H^{0}(\cdot)=\emptyset$).\\

Let us now specify the regions on which depend the growing segments. Let $x=(\xi,u)\in\mathbb{R}^{2}\times [0;1]$ and $r>0$. Here, the crucial point is to remark that the indicator function
\begin{equation}
\label{IndicatorLocal}
\1_{\{\|\xi-h_{g}(\mathbb{X}\cup\{x\},x)\| \le r\}} \, \mbox{ is $\mathscr{S}_{B(\xi+r\overrightarrow{u},r)}$-measurable.}
\end{equation}
See Figure \ref{fig:stability}. Let $\varphi\in\mathscr{C}'$ and $\Lambda$ be an open bounded region in $\mathbb{R}^{2}$. For each marked point $x=(\xi,u)\in\varphi_{\Lambda}$, we set
$$
r(x,\Lambda) = \sup \{ r\ge 0 , \, B(\xi+r\overrightarrow{u},r) \subset \Lambda \} ~.
$$
Hence, for all $r\le r(x,\Lambda)$, it is possible to know when we observe the configuration $\varphi$ only through the window $\Lambda$ if $\|\xi-h_{g}(\varphi,x)\|$ is smaller than $r$ or not. Henceforth, we define the \textit{decision set} of the growing segment $x$ through $\Lambda$ as
$$
D_{\Lambda}(x) = B(\xi+r(x,\Lambda)\overrightarrow{u},r(x,\Lambda)) ~.
$$

\begin{figure}[!ht]
\begin{center}
\psfrag{x}{\small{$x$}}
\psfrag{y1}{\small{$y_{1}$}}
\psfrag{y2}{\small{$y_{2}$}}
\psfrag{Lambda}{\small{$\Lambda$}}
\psfrag{Dlambdax}{\color{red}{\small{$D_{\Lambda}(x)$}}}
\includegraphics[width=5.5cm,height=5cm]{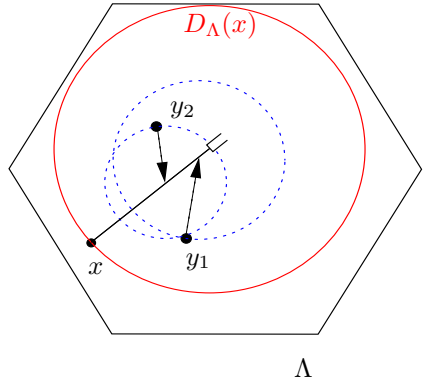}
\end{center}
\caption{  \label{fig:stability}  The set $\Lambda$ is a regular hexagon. The decision set $D_{\Lambda}(x)$  is delimited by a red circle. It contains two marked points $y_{1},y_{2} \in \varphi$  which could stop the growing segment $x$ before $\xi+r(x,\Lambda) u  $. To check it, it is enough to observe the configuration $\varphi$ inside the blue dashed circles which are themselves included in $D_{\Lambda}(x)$.}
\end{figure}

So, for a given marked point $x=(\xi,u)\in\varphi_{\Lambda}$, two situations may occur. If the stopping vertex of $x$ in $\varphi$ belongs to the decision set $D_{\Lambda}(x)$ then the whole segment $[\xi,h_{g}(\varphi,x)]$ is observed knowing $\varphi_{\Lambda}$. In this case, we set
$f_{\Lambda}(\varphi,x)=h_{g}(\varphi,x)$. Otherwise, we can only assert that the line-segment $x$ will be longer than $r(x,\Lambda)$. In that case, $f_{\Lambda}(\varphi,x)=\xi+r(x,\Lambda)\overrightarrow{u}$. In both situations,
$$
\left[ \xi;f_{\Lambda}(\varphi,x) \right] \subset \left[ \xi;h_{g}(\varphi,x) \right] ~.
$$
The previous considerations lead to the next result:

\begin{lem}
\label{lem:LocalDep}
With the above notations, the random set
$$
\mathscr{G}_{\Lambda}(\varphi) := \bigcup_{x\in\varphi_{\Lambda}} \left[ \xi;f_{\Lambda}(\varphi,x) \right]
$$
is $\mathscr{S}_{\Lambda}$-measurable.
\end{lem}

\begin{dem}
Let $x=(\xi,u)\in\varphi_{\Lambda}$. By construction, the random segment $[\xi;f_{\Lambda}(\varphi,x)]$ contains all the information available on the segment $x=(\xi,u)$ when we only observe $\varphi$ through the decision set $D_{\Lambda}(x)=B(\xi+r(x,\Lambda)\overrightarrow{u},r(x,\Lambda))$: either it is still alive at $\xi+r(x,\Lambda)\overrightarrow{u}$ or it has been stopped before. By definition of $r(x,\Lambda)$, $[\xi;f_{\Lambda}(\varphi,x)]$ is $\mathscr{S}_{\Lambda}$-measurable.
\end{dem}

We can now introduce the central notion of shield hexagons.The hexagon $\H(\xi)$ is $\epsilon$-shield (for $\varphi$) whenever the set $\mathscr{G}_{\H(\xi)}(\varphi)$ produces a barrier in the strip $\H(\xi)\setminus (\xi+\epsilon\H(0))$ disconnecting the inside part $\xi+\epsilon\H(0)$ from the outside part $\H(\xi)^{c}$. 

\begin{defin}
Let $\varphi\in\mathscr{C}^{'}$, $\epsilon\in (0,1)$ and $\xi\in\mathbb{R}^{2}$. The hexagon $\H(\xi)$ is said \textrm{$\epsilon$-shield} for $\varphi$ if for all $a,b\in\mathbb{R}^{2}$ such that $a\notin\H(\xi)$ and $b\in \xi+\epsilon\H(0)$, we have
$$
(a;b) \cap \mathscr{G}_{\H(\xi)}(\varphi) \not= \emptyset ~.
$$
Moreover, for any integer $n>0$ and $\{z_{i}\}_{1\le i\le n}\subset\Pi$, the collection $\{\H(z_{i})\}_{1\le i\le n}$ is said \textrm{$\epsilon$-shield} for $\varphi$ if for each index $i$, $\H(z_{i})$ is $\epsilon$-shield for $\varphi$.
\end{defin}

It is not difficult to be convinced (using many small segments, all the smaller as $\epsilon\to 1$) that this event occurs with positive probability:
\begin{equation}
\label{ProbaHexShield>0}
\forall \epsilon \in (0;1) , \, p_{\epsilon} = \mathbb{P} \left[ \H(0) \text{ is $\epsilon$-shield} \right] \, > \, 0 ~.
\end{equation}

The notion of decision sets $D_{\Lambda}(\cdot)$-- and also $\mathscr{G}_{\Lambda}(\cdot)$ --have been introduced to use the independence property of the Poisson point process $\mathbb{X}$. Indeed, by Lemma \ref{lem:LocalDep}, for any vertices $z\not=z^{'}\in\Pi$, the hexagons $\H(z)$ and $\H(z')$ are independently $\epsilon$-shield.\\

The next step consists in using $\epsilon$-shield hexagons as building blocks to create obstacles. Precisely:

\begin{defin}
\label{defin:Shielded}
Let $m\in\mathbb{N}^{*}$ be an integer and $\varphi\in\mathscr{C}^{'}$ a marked configuration. Any $\eta\in\mathbb{R}^{2}$ is said $m$-\textrm{shielded} for $\varphi$ if:
\begin{itemize}
\item[$(\clubsuit)$]For all $x=(\xi,u)\in\varphi_{\H^{2m}(\eta)^{c}}$, $[\xi;h_{g}(\varphi,x)]\cap \H^{m}(\eta)=\emptyset$;
\item[$(\spadesuit)$]For all $x\in\varphi_{\H^{m}(\eta)}$, $h_{g}(\varphi,x)\in\H^{2m}(\eta)$.
\end{itemize}
\end{defin}

If $\eta$ is $m$-shielded for $\varphi$ then conditions $(\clubsuit)$ and $(\spadesuit)$ roughly asserts that it is impossible for a growing segment to cross $\H^{2m}(\eta)\setminus \H^{m}(\eta)$ respectively from the outside part $\H^{2m}(\eta)^{c}$ and from the inside part $\H^{m}(\eta)$.

In the sequel, we will establish the existence of an event $E_{m}\in\mathscr{S}_{Hex^{2m}(0)}$ such that, on $E_{m}$, $0$ is a.s. $m$-shielded (Proposition \ref{prop:OShielded}). Actually, it will be required to get the event $E_{m}$ that $\H^{2m}(0)\setminus\H^{m}(0)$ contains many $\epsilon$-shield hexagons. In a second time, we will prove that the probability of $E_{m}$ tends to $1$ as $m\to\infty$ (Proposition \ref{prop:LimPrEm3m}).

Let us introduce some notations needed to define the event $E_{m}$. Let $\eta\in\partial\H^{m}(0)$ where $\partial\Lambda$ denotes the topological boundary of $\Lambda\subset\mathbb{R}^{2}$. For any $v\in[0;1]$, we define the (semi-infinite) ray starting from $\eta$ in the direction $\overrightarrow{v}=(\cos(2\pi v),\sin(2\pi v))$ by $l(\eta,\overrightarrow{v})=\{\eta+t\overrightarrow{v}, t\geq 0\}$. Thus, we denote by $\mathscr{L}^{m}$ the set of rays $l(\eta,\overrightarrow{v})$ coming from $\partial\H^{m}(0)$ which do not overlap the topological interior of $\H^{m}(0)$:
$$
\mathscr{L}^{m} = \left\{ l(\eta,\overrightarrow{v}) ,\ l(\eta,\overrightarrow{v}) \cap \text{Int}(\H^{m}(0)) = \emptyset \, \mbox{ and } \, (\eta,v)\in\partial\H^{m}(0)\times[0;1] \right\} ~.
$$
For each ray $l\in\mathscr{L}^{m}$, let us consider the set of hexagons included in $\H^{2m}(0)\setminus\H^{m}(0)$ and crossed by $l$:
$$
\text{Cross}(l) = \left\{ \H(z) ,\ m+1\le d_{\Pi}(0,z)\le 2m \, \mbox{ and } \, l\cap\H(z) \not= \emptyset \right\} ~.
$$
This set can be partitioned into different floors $\text{Cross}_{i}(l)$, for $m+1\leq i\leq 2m$, where $\text{Cross}_{i}(l)$ denotes the set of hexagons of $\text{Cross}(l)$ included in $C_{i}(0)$.
We can observe that, for each $l\in\mathscr{L}^{m}$, there exists $m+1\le i(l)\le 2m$ such that for all $i(l)\le i \le 2m$, $\text{Cross}_{i}(l)$ contains at most three hexagons.\\
Thus, the set $\text{Cross}(l)$ is said $\epsilon$\textit{-uncrossable} for $\varphi$ if we can find two consecutive floors $\text{Cross}_{i}(l)$ and $\text{Cross}_{i+1}(l)$, for some index $i(l)\leq i\leq 2m-1$, which are both $\epsilon$-shield for $\varphi$. We can then define the event $E_{m}$ as:
\begin{equation}
\label{Emn}
E_{m}(\epsilon) = \bigcap_{l\in\mathscr{L}^{m}} \left\{ \text{Cross}(l)\text{ is }\epsilon\text{-uncrossable for }\mathbb{X} \right\}
\end{equation}
which is $\mathscr{S}_{Hex^{2m}(0)}$-measurable by construction.

\begin{prop}
\label{prop:OShielded}
There exists $\epsilon\in (0,1)$ (close to $1$) such that, a.s. on the event $E_{m}(\epsilon)$, $0$ is $m$-shielded.
\end{prop}

\begin{dem}
Assume that the event $E_{m}(\epsilon)$ is satisfied. Then, to prove that $0$ is $m$-shielded, i.e. to prove $(\clubsuit)$ and $(\spadesuit)$, it is enough to state that any ray $l\in\mathscr{L}^{m}$ crosses the inside part $z+\epsilon\H(0)$ of an $\epsilon$-shield hexagon $\H(z)$.

Let us consider a ray $l\in\mathscr{L}^{m}$ and a vertex $z\in B_{2m}(0)\!\setminus\!B_{i(l)-1}(0)$. Let us define $d_{z,l}$ as follows:
$$
d_{z,l} = \sup_{x\in l\cap\H(z)} d(x, \partial\H(z))
$$
and $d_{z,l}=0$ if $l\cap\H(z)$ is empty (where the above distance $d$ is euclidean). Thus, we set
$$
\gamma = \inf_{l\in\mathscr{L}^{m}} \sup \{ d_{z,l}\, , \ z\in B_{n}(0)\!\setminus\!B_{i(l)-1}(0) \, \mbox{ and $\H(z)$ is $\epsilon$-shield} \} ~.
$$
On the event $E_{m}(\epsilon)$, any ray $l$ crosses two consecutive $\epsilon$-shield floors. Hence, the hexagonal construction ensures that the above infimum $\gamma$ is positive (see Figure \ref{fig:hexagoncross}). So, $\epsilon=1-\gamma/2$ is suitable.
\end{dem}

\begin{figure}[!ht]
\begin{center}
\psfrag{z}{\small{$z$}}
\psfrag{z'}{\small{$z'$}}
\psfrag{l}{\small{$l$}}
\psfrag{dm}{\small{$\partial \H^{m}(0)$}}
\psfrag{d3m}{\small{$\partial \H^{2m}(0)$}}
\includegraphics[width=10cm,height=6cm]{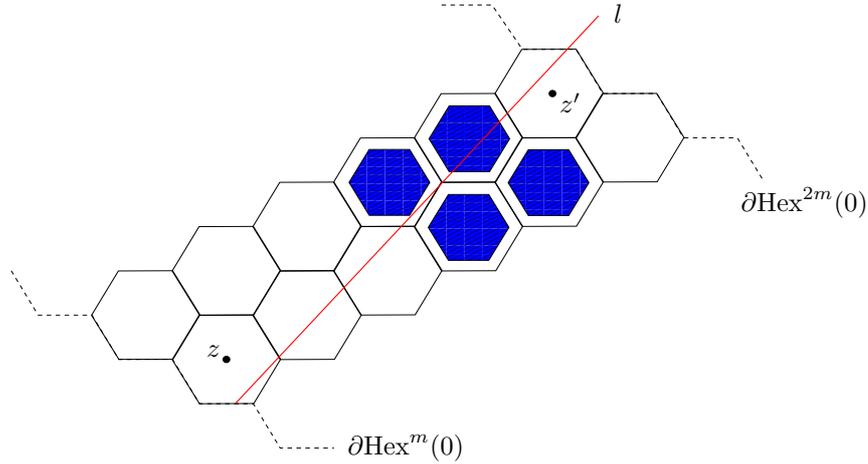}
\end{center}
\caption{\label{fig:hexagoncross} Here is a ray $l\in\mathscr{L}^{m}$ (in red) starting at $\partial \H^{m}(0)$ and crossing two consecutive $\epsilon$-shield floors, say $\text{Cross}_{i}(l)$ and $\text{Cross}_{i+1}(l)$. Each hexagon $\H(z)$ belonging to these two floors has its inside part, i.e. $z+\epsilon\H(0)$, colored in blue. As this picture suggests, when $\epsilon$ is close to $1$, it becomes impossible for the ray $l$ to avoid the inside parts of hexagons of two consecutive floors.}
\end{figure}
In the sequel, we merely write $E_{m}$ instead of $E_{m}(\epsilon)$ where $\epsilon$ is given by Proposition \ref{prop:OShielded}. Its probability tends to $1$ with $m$;

\begin{prop}
\label{prop:LimPrEm3m}
The probability of the event $E_{m}$ tends to 1 as $m$ tends to $+\infty$.
\end{prop}

\begin{dem}
First, let us reduce the infinite intersection defining the event $E_{m}$ in (\ref{Emn}) to a finite one. Let $z\in S^{m+1}(0)$ and $z'\in S^{2m}(0)$. Let us consider the set $\text{Cross}(z,z')$ made up of hexagons $\H(z'')$, $z''\in\Pi$, which are crossed by a ray $l\in\mathscr{L}^{m}$ starting at some $\eta\in\partial\H(z)$ and exiting $\H^{2m}(0)$ through $\H(z')$. As previously, we divide the set $\text{Cross}(z,z')$  into different floors $\text{Cross}_{i}(z,z')$, for $m+1\leq i\leq 2m$, where $\text{Cross}_{i}(z,z')$ denotes the set of hexagons of $\text{Cross}(z,z')$ included in $C_{i}(0)$.\\
There exists $m+1\le i(z,z')\le 2m$ such that for all $i(z,z')\le i \le 2m$, $\text{Cross}_{i}(z,z')$ contains at most three hexagons.Thus, $\text{Cross}(z,z')$ is said $\epsilon$\textit{-uncrossable} for $\mathbb{X}$ if we can find two consecutive floors $\text{Cross}_{i}(z,z')$ and $\text{Cross}_{i+1}(z,z')$, for some index $i(l)\leq i\leq 2m-1$, which are both $\epsilon$-shield for $\mathbb{X}$.
 Hence,
\begin{equation}
\label{inclusionEm3m}
\bigcap_{z\in S^{m+1}(0), z'\in S^{2m}(0)} \left\{ \text{Cross}(z,z')\text{ is }\epsilon\text{-uncrossable for }\mathbb{X} \right\} \, \subset \, E_{m} ~.
\end{equation}
Since $S^{m+1}(0)\times S^{2m}(0)$ contains $36(m+1)(2m)$ vertices, the expected result follows from the inclusion (\ref{inclusionEm3m}) and the limit: for $m$ sufficiently large, for all $z\in S^{m+1}(0)$ and $z'\in S^{2m}(0)$,
\begin{equation}
\label{LimCrosszz'}
\mathbb{P}\left[ \text{Cross}(z,z')\text{ is not $\epsilon$-uncrossable for $\mathbb{X}$} \right] \leq (1- p_{\epsilon}^{6})^{\frac{m}{10}} ~,
\end{equation}
where $p_{\epsilon}=\mathbb{P}[\H(0) \text{ is $\epsilon$-shield}]$.

We obtain
$$\lbrace \text{Cross}(z,z^{'})\text{ is not }\epsilon-\text{uncrossable for }\mathbb{X}\rbrace\subset\bigcap_{k=i(z,z')}^{2m-1} U_{k}$$
where $U_{k}=\lbrace \text{Cross}_{k}(z,z')$ and $\text{Cross}_{k+1}(z,z')\text{ are }\epsilon-\text{ shield for }\mathbb{X}\rbrace ^{c}.$\\
To obtain an independence property under the Poisson point process law, we need to consider disjoint subsets of hexagons:
$$T_{m}=\bigcap_{k=0}^{\lfloor\frac{2m-1-i(z,z')}{2}\rfloor}U_{i(z,z')+2k} .$$
 The events $(U_{2k})_{k}$ are mutually independent and  
$$\lbrace \text{Cross}(z,z^{'})\text{ is not }\epsilon-\text{uncrossable for } \mathbb{X}\rbrace \subset T_{m}.$$
We have introduced in (\ref{ProbaHexShield>0}) the probability $p_{\epsilon}=\mathbb{P}\left[ \H(0)\text{ is }\epsilon-\text{shield for }\mathbb{X}\right]$. Then, for $i(z,z')\le k\le 2m$, we have $\mathbb{P}[U_{k}]\le 1- p_{\epsilon}^{6}$. It is relatively easy to check that, for $m$ sufficiently large, for all $(z,z')\in S^{m+1}(0)\times S^{2m}(0)$, we have: $\lfloor\frac{2m-1-i(z,z')}{2}\rfloor +1\ge\frac{m}{10} $. It implies the existence of a bound for $\mathbb{P}[T_{m}]$:
  $$\mathbb{P}\left[T_{m}\right]\le ( 1- p_{\epsilon}^{6})^{\frac{m}{10}}.$$
 Then,
 $$\mathbb{P}\left[ \text{Cross}(z,z^{'})\text{ is not }\epsilon-\text{shield for } \mathbb{X}\right]\le ( 1- p_{\epsilon}^{6})^{\frac{m}{10}}.$$
\end{dem}

From now on, we claim that:

\begin{prop}
The line-segment model satisfies \textbf{(SA)} for $\alpha=32$ and $\mathscr{E}_{m}=E_{m}\cap E_{2m}$.
\end{prop}

\begin{dem}
We have to check that the line-segment model satisfies the three items of \textbf{(SA)}. By definition, the event $\mathscr{E}_{m}=E_{m}\cap E_{2m}\in\mathscr{S}_{\H^{4m}(0)}$. Item $(i)$ follows from the fact that any $\eta\in\H^{4m}(0)$ satisfies $\|\eta\|\leq 4m\sqrt{3}+1\leq 8m$. Item $(ii)$ is given by Proposition \ref{prop:LimPrEm3m}. So, it only remains to check Item $(iii)$.

For this purpose, let us consider three disjoint subsets $V,A_{1},A_{2}$ of $\mathbb{Z}^{2}$ such that $\partial A_{1}$ and $\partial A_{2}$ are included in $V$. Let also for $i\in\{1,2\}$,
$$
\mathcal{A}_{i} = \left( A_{i}\oplus[-\frac{1}{2},\frac{1}{2}]^{2} \right) \setminus (V\oplus[-\alpha,\alpha]^{2}) ~.
$$
Thus, let $m$ be a positive integer and $\varphi,\varphi'\in\mathscr{C}'$ such that $\tau_{-mz}(\varphi)\in\mathscr{E}_{m}$, for all $z\in V$. We have to check that
\begin{equation}
\label{item3}
\forall x \in \varphi_{m\mathcal{A}_{1}} , \,  h(\varphi,x) = h(\bar{\varphi},x) ~,
\end{equation}
where $\bar{\varphi}$ denotes the configuration $\varphi_{m\mathcal{A}_{2}^{c}}\cup \varphi'_{m\mathcal{A}_{2}}$. The reason why (\ref{item3}) holds can be roughly expressed as follows. The replacement of the configuration $\varphi$ with $\bar{\varphi}$, which actually concerns only the set $m\mathcal{A}_{2}$, may generate some modifications in the graph on the set $m\mathcal{A}_{2}^{c}$ but not beyond the obstacle $mV$. Then, the graph on $m\mathcal{A}_{1}$ is preserved.

Let us start with splitting the set $m\mathcal{A}_{2}^{c}$ into three disjoint subsets: $m\mathcal{A}_{1}$, $\text{Shield}:=mV\oplus\H^{2m}(0)$ and $\text{Bound}:=\big( mV\oplus [-\alpha m;\alpha m]^{2}\big)\setminus\text{Shield}$. Since $mV\oplus\H^{4m}(0)$ is included in $\text{Bound}$ and $\mathscr{E}_{m}$ is $\mathscr{S}_{\H^{4m}(0)}$-measurable then
\begin{equation}
\label{StillEm}
\forall z \in V , \, \tau_{-mz}(\bar{\varphi}) \in \mathscr{E}_{m} ~.
\end{equation}
In other words, the shield structure of vertices $mz$, $z\in V$, is preserved when passing from $\varphi$ to $\bar{\varphi}$. Hence, any $x=(\xi,\cdot)\in\varphi_{\text{Shield}}$ belongs to a set $\H^{2m}(mz)$ where $mz\in mV$ is $2m$-shielded (thanks to $E_{2m}$). Property $(\spadesuit)$ of Definition \ref{defin:Shielded} then ensures that $h_{g}(\varphi,x)$ is in $mz\oplus\H^{4m}(0)$, i.e. $\|\xi-h_{g}(\varphi,x)\|\le \text{diam}\, \H^{4m}(0)\leq 16m$. To sum up, for all $x\in\varphi_{\text{Shield}}$,
$$
B(h_{g}(\varphi,x),\|\xi-h_{g}(\varphi,x)\|) \subset mV\oplus [-32m;32m]^{2} ~.
$$
The above inclusion justifies the choice $\alpha=32$. It also ensures that the replacement of $\varphi$ with $\varphi'$ outside of $mV\oplus [-32m;32m]^{2}$ does not impact the geometric edges starting from vertices of $\varphi_{\text{Shield}}$. So,
\begin{equation}
\label{Shieldnotcross}
\forall x \in \varphi_{\text{Shield}} , \, h(\bar{\varphi},x)= h(\varphi,x) ~.
\end{equation}

It then remains to show that (\ref{item3}) is a consequence of (\ref{StillEm}) and (\ref{Shieldnotcross}). When passing from $\varphi$ to $\bar{\varphi}$, the line-segment of a marked point $x\in\varphi_{m\mathcal{A}_{2}^{c}}$ can be modified in two different ways:
\begin{itemize}
\item[$\bullet$] Either the line-segment of $x$ is shorter for $\bar{\varphi}$ than for $\varphi$, i.e. $x$ admits a new outgoing neighbor $y$.
\item[$\bullet$] Or the line-segment of $x$ is longer for $\bar{\varphi}$ than for $\varphi$, i.e. its original stopping line-segment has been stopped before by some marked point $y$.
\end{itemize}
In both cases, we say that the marked point $y$ modifies $x$. It belongs to $\varphi'_{m\mathcal{A}_{2}}$, or to $\varphi_{m\mathcal{A}_{2}^{c}}$ but with the condition that $h(\varphi,y)\not= h(\bar{\varphi},y)$. Hence, from any $x_{0}\in\varphi'_{m\mathcal{A}_{2}}$, may start a sequence of marked points $(x_{i})_{0\leq i\leq n}$ such that $x_{i}$ modifies $x_{i+1}$. Now, (\ref{StillEm}) and (\ref{Shieldnotcross}) prevents such sequence to cross the set $\text{Shield}=mV\oplus\H^{2m}(0)$. By contradiction, let us assume that $x_{n}\in m\mathcal{A}_{1}$. Since (\ref{Shieldnotcross}) prevents the $x_{i}$'s to belong to $\text{Shield}$, it necessarily exists an index $0\leq i_{0}<n$ such that the line-segment of $x_{i_{0}}$ crosses $mV\oplus\H^{m}(0)$. But this is forbidden by (\ref{StillEm}): each $mz$, for $z\in V$, is $m$-shielded (thanks to $E_{m}$). So Property $(\clubsuit)$ of Definition \ref{defin:Shielded} applies.
\end{dem}

\subsection{Navigation model}

Let $0<\epsilon<\frac{\pi}{2}$. Given a configuration $\varphi\in\mathscr{C}^{'}$, let us recall that the stopping vertex of $x=(\xi,u)\in\varphi$ is the closest element of $\varphi_{\text{germs}}\cap C(x)$ to $\xi$, where
$$
C(x)=\left\{ (r\cos(\alpha),r\sin(\alpha)) ; \, r>0 \; \mbox{ and } \; |\alpha - 2\pi u| < \epsilon \right\} ~.
$$
If $(\eta,v)=h(\varphi,x)$ then the impact point of $x$ in the Navigation model is $h_{g}(\varphi,x)=\eta$.

\subsubsection{Loop assumption}

The Navigation model satisfies \textbf{(LA)}.

\begin{prop}
\label{Prop:CycAssNav}
Each configuration of $\mathscr{C}^{'}$ is $1-\text{looping}$.
\end{prop}

\begin{dem}
Let $\varphi\in\mathscr{C}^{'}$ and $x=(\xi,u)\in\varphi$. Let us introduce the stopped cone starting from $x$:
$$
C_{\text{stop}}(x) = \left\{ (r\cos(\alpha),r\sin(\alpha)) ; \, 0<r<\|\xi-h_{g}(\varphi,x)\| \,\ \mbox{ and } \; |\alpha - 2\pi u| < \epsilon \right\} ~.
$$
Hence, $C_{\text{stop}}(x)\cap\varphi_{\text{germs}}=\emptyset$. Let $d_{x}>0$ small enough so that $\varphi_{B(\xi,d_{x})}$ only contains the point $\xi$. Thus, let us consider an open ball $A_{x}\subset \mathbb{R}^{2}\times[0;1]$ such that any marked point $y=(\eta,v)\in A_{x}$ satisfies $\eta\in B(\xi,d_{x})\cap C_{stop}(x)$ and $\xi\in C(y)$. Thenceforth, we get $h(\varphi\cup\{y\},x)=y$ and $h(\varphi\cup\{y\},y)=x$ (see Figure \ref{fig:cyclenav}) which respectively imply
$$
\text{For}(x,\varphi\cup\{y\}) = \{x,y\} \; \mbox{ and } \; \text{Back}(x,\varphi)\cup\{y\} \subset \text{Back}(x,\varphi\cup\{y\}) ~.
$$
Let us justify this latter inclusion. It is possible that $\eta$ belongs to the stopped cone $C_{\text{stop}}(z)$ (w.r.t. $\varphi$) of a given marked point $z\in\varphi$, which forces $h(\varphi\cup\{y\},z)=y$. Since $y$ belongs to $\text{Back}(x,\varphi\cup\{y\})$, the same holds for $z$. If $z$ was already in the backward of $x$ (for $\varphi$), it is still in (but for $\varphi\cup\{y\}$). See Figure \ref{fig:cyclenav}.
\end{dem}

\begin{figure}[!ht]
\begin{center}
\psfrag{x}{\small{$x$}}
\psfrag{h}{\small{$h(\varphi,x)$}}
\psfrag{z}{\small{$z$}}
\psfrag{y}{\small{\color{red}{$y$}}}
\includegraphics[width=8.5cm,height=5cm]{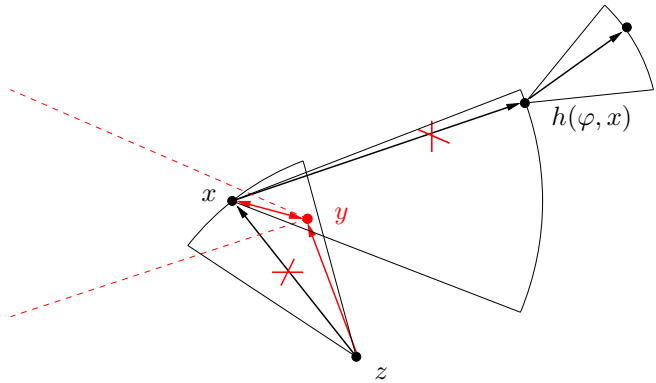}
\end{center}
\caption{\label{fig:cyclenav} The adding of the marked point $y$ breaks the edges from $x$ to $h(\varphi,x)$, and from $z$ to $x$. However, for the configuration $\varphi\cup\{y\}$, $z$ is still in the backward of $x$.}
\end{figure}

\subsubsection{Shield assumption}

Let us split the square $[-m;m]^{2}$ into $\kappa=(2\lfloor m^{1/2}\rfloor)^{2}$ congruent subsquares $Q^{m}_{1},\ldots,Q^{m}_{\kappa}$ ($\lfloor\cdot\rfloor$ denotes the integer part). Each of these subsquares has an area equal to
$$
\left( \frac{2m}{2\lfloor m^{1/2}\rfloor} \right)^{2} ~,
$$
i.e. of order $m$. Thus, we define the event $\mathscr{E}_{m}$ as follows:
$$
\mathscr{E}_{m} = \bigcap_{1\leq i\leq\kappa} \left\{ \#\mathbb{X}_{Q^{m}_{i}}\ge 1 \right\} ~.
$$

\begin{prop}
\label{Prop:ShiAssNav}
For $\alpha=1$, the Navigation model satisfies \textbf{(SA)} w.r.t. the  family of events $(\mathscr{E}_{m})_{m\ge 1}$.
\end{prop}

\begin{figure}[!ht]
\begin{center}
\psfrag{x}{\small{$x=(\xi,u)$}}
\includegraphics[width=8cm,height=5.5cm]{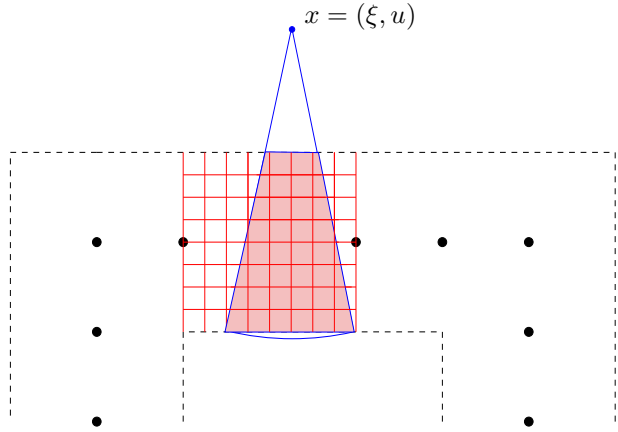}
\end{center}
\caption{\label{fig:shieldnav} Black points are vertices $mz$ for $z\in V$. The event $\mathscr{E}_{m}$ realized on each $z\oplus [-m;m]^{2}$ provides a shield between $m\mathcal{A}_{1}$ and $m\mathcal{A}_{2}$. Indeed, the cone $C_{\text{stop}}(x)$ cannot overlap $m\mathcal{A}_{2}$ without containing a subsquare $z+Q_{i}^{m}$.}
\end{figure}

\begin{dem}
Let us first remark that the event $\mathscr{E}_{m}$ is $\mathscr{S}_{[-m,m]^{2}}$-measurable and its probability tends to $1$. So the first two items of \textbf{(SA)} are satisfied with $\alpha=1$.

Let us focus on Item (iii). Hence, let us consider $V,A_{1},A_{2}\subset\mathbb{Z}^{2}$ such that the topological conditions of \textbf{(SA)} occur. Let us set
$$
\mathcal{A}_{i} = \left( A_{i}\oplus [-\frac{1}{2},\frac{1}{2}]^{2} \right) \setminus (V\oplus[-\alpha,\alpha]^{2}) 
$$
for $i\in\{1,2\}$. Let $\varphi\in\mathscr{C}^{'}$ satisfying $\varphi-mz\in\mathscr{E}_{m}$, for any vertex $z\in V$. Let $x=(\xi,u)\in\varphi$ be a marked point whose first coordinate belongs to $m\mathcal{A}_{i}$. If the cone $C(x)$ does not overlap $m\mathcal{A}_{j}$, with $j=3-i$, then the outgoing vertex $h(\varphi,x)$ does not depend on possible changes on $\varphi_{m\mathcal{A}_{j}}$. From now on, let us assume that $C(x)\cap m\mathcal{A}_{j}$ is not empty (see Figure \ref{fig:shieldnav}). It is then sufficient to remark that for any $m\geq m_{0}(\epsilon)$, the stopped cone $C_{\text{stop}}(x)$ does not overlap $m\mathcal{A}_{j}$. Otherwise, for $m$ large enough, it would contain at least one subsquare $z+Q_{i}^{m}$ for some $z\in V$ and $1\leq i\leq\kappa$ and so at least a marked point (since $\varphi-mz\in\mathscr{E}_{m}$) which is forbidden. Hence, as previously, the outgoing vertex $h(\varphi,x)$ remains unchanged whatever the configuration $\varphi$ inside $m\mathcal{A}_{j}$. 
\end{dem}

\subsubsection{Acknowledgement}

The authors thank C. Hirsch, G. Last and S. Zuyev for fruitful discussions on this problem. {They also thank J.-B. Gou\'er\'e for bringing to our attention an example of outdegree-one graph which satisfies \textbf{(LA)}, but not \textbf{(SA)}, and which percolates.}This work was supported in part by the Labex CEMPI (ANR-11-LABX-0007-01), the CNRS GdR 3477 GeoSto and the ANR PPPP (ANR-16-CE40-0016).

\addcontentsline{toc}{section}{Bibliography}

\bibliographystyle{plain}
\bibliography{Biblio}
\end{document}